\documentclass[a4paper,12pt, twoside,leqno]{article}

\usepackage[margin=25mm]{geometry}

\usepackage[utf8]{inputenc}
\usepackage[T1]{fontenc}
\usepackage[british]{babel}
\usepackage[style=british]{csquotes}
\usepackage[en-GB]{datetime2}
\DTMlangsetup[en-GB]{ord=omit}
\usepackage{xcolor}

\usepackage{amsthm, thmtools}
\usepackage{mathtools}
\usepackage{titlesec} 
\usepackage{aliascnt}
\usepackage{hyperref}
\hypersetup{
    bookmarksnumbered,
    colorlinks,
    linkcolor={purple!50!black},
    citecolor={blue!50!black},
    urlcolor={green!50!black}
}
\usepackage[capitalise,nameinlink]{cleveref}

\usepackage{tikz}
\usepackage{tikz-cd}

\usepackage[fulladjust]{marginnote}
\usepackage{stmaryrd} 
\usepackage[textsize=scriptsize]{todonotes}

\usepackage{aliascnt}

\usepackage{array}

\usepackage{quiver}

\usepackage{newtxtext}
\usepackage[lining,tabular,scale=.95]{FiraSans}
\usepackage{newtxmath}
\usepackage[lining,scaled=.8]{FiraMono} 
\usepackage[cal=boondoxo,frak=boondox]{mathalpha} 

\newcommand*{\titlefont}{\sffamily\scshape}
\newcommand*{\subtitlefont}{\firamedium}

\DeclareSymbolFont{oldtxsymbols}{OMX}{txex}{m}{n}
\DeclareMathSymbol{\intop}{\mathop}{oldtxsymbols}{"52}
\let\int\intop

\DeclareFontEncoding{MDA}{}{}
\DeclareSymbolFont{mathdesignA}{MDA}{mdput}{b}{n} 
\SetSymbolFont{mathdesignA}{bold}{MDA}{mdput}{b}{n}
\DeclareSymbolFontAlphabet{\mathbb}{mathdesignA}

\tikzcdset{
    arrow style=tikz,
    arrows={/tikz/line width=.5pt},
    diagrams={>={Straight Barb[scale=0.8]}}
}

\DeclareFontFamily{OMX}{MnSymbolE}{}
\DeclareSymbolFont{MnLargeSymbols}{OMX}{MnSymbolE}{m}{n}
\SetSymbolFont{MnLargeSymbols}{bold}{OMX}{MnSymbolE}{b}{n}
\DeclareFontShape{OMX}{MnSymbolE}{m}{n}{
    <-6>  MnSymbolE5
   <6-7>  MnSymbolE6
   <7-8>  MnSymbolE7
   <8-9>  MnSymbolE8
   <9-10> MnSymbolE9
  <10-12> MnSymbolE10
  <12->   MnSymbolE12
}{}
\DeclareFontShape{OMX}{MnSymbolE}{b}{n}{
    <-6>  MnSymbolE-Bold5
   <6-7>  MnSymbolE-Bold6
   <7-8>  MnSymbolE-Bold7
   <8-9>  MnSymbolE-Bold8
   <9-10> MnSymbolE-Bold9
  <10-12> MnSymbolE-Bold10
  <12->   MnSymbolE-Bold12
}{}
\DeclareMathDelimiter{[}{\mathopen}{MnLargeSymbols}{'000}{MnLargeSymbols}{'000}
\DeclareMathDelimiter{]}{\mathclose}{MnLargeSymbols}{'005}{MnLargeSymbols}{'005}
\DeclareMathDelimiter{\llbr}{\mathopen}{MnLargeSymbols}{'102}{MnLargeSymbols}{'102}
\DeclareMathDelimiter{\rrbr}{\mathclose}{MnLargeSymbols}{'107}{MnLargeSymbols}{'107}

\usepackage{bm}
\usepackage{cases}
\usepackage{accents}

\usepackage{setspace}
\usepackage{subdepth}

\newcommand{\initlengths}{%
    \setlength{\abovedisplayshortskip}{3pt plus 9pt minus 3pt}%
    \setlength{\belowdisplayshortskip}{9pt plus 9pt minus 9pt}%
    \setlength{\abovedisplayskip}{9pt plus 9pt minus 9pt}%
    \setlength{\belowdisplayskip}{9pt plus 9pt minus 9pt}%
    \tolerance 500
}

\usepackage{titling} 
\pretitle{\vspace{-30pt}\setstretch{1.05}\begin{center}\LARGE\scshape\sffamily\fontdimen2\font=0.3em} 
\posttitle{\end{center}\vspace{-3pt}} 

\numberwithin{paragraph}{subsection} 
\newcommand{\parafont}{\bfseries} 
\newcommand{\parasep}{9pt plus 3pt minus 3pt} 
\titleformat{\section}{\center\Large\titlefont}{\thesection}{1em}{} 
\titleformat{\subsection}{\large\subtitlefont\boldmath}{\thesubsection}{1em}{} 
\titleformat{\paragraph}[runin]{\parafont}{\theparagraph.}{.33em}{\normalfont\bfseries\boldmath} 
\titlespacing{\paragraph}{0pt}{\parasep}{.5em} 

\usepackage{tocloft} 
\renewenvironment{abstract}{%
    \centering\begin{minipage}{.85\textwidth}%
    \setlength{\parindent}{1.5em}%
    \centerline{\subtitlefont\abstractname}%
    \setstretch{1.15}%
    \par\vspace{6pt}%
}{\end{minipage}\par\vspace{3pt}}

\newcommand{\authorinforule}{\noindent\rule{0.38\textwidth}{0.4pt}}

\newlength{\authorwidth}
\newcommand{\authorinfo}[3]{%
    \setlength{\leftskip}{1.5em}
    \setlength{\parindent}{0em}
    \par%
    {\small%
    \makebox[\authorwidth][l]{#1}%
    \texttt{#2}%
    \\[-2pt]
    #3.}
    \vspace{4pt}\par
}

\makeatletter
\declaretheoremstyle[
    spaceabove=\parasep, spacebelow=\parasep,
    postheadspace=.5em,
    headfont=\normalfont\bfseries,
    headpunct={},
    headformat={\NUMBER.\@\ \NAME.\@\NOTE},
    notefont=\normalfont\bfseries\boldmath,
    notebraces={}{.},
    bodyfont=\itshape,
]{theorem}
\declaretheoremstyle[
    spaceabove=\parasep, spacebelow=\parasep,
    postheadspace=.5em,
    headfont=\normalfont\bfseries,
    headpunct={},
    headformat={\NAME.\@\NOTE},
    notefont=\normalfont\bfseries\boldmath,
    notebraces={}{.},
    bodyfont=\itshape,
]{theorem*}
\declaretheoremstyle[
    spaceabove=\parasep, spacebelow=\parasep,
    postheadspace=.5em,
    headfont=\normalfont\bfseries,
    headpunct={},
    headformat={\NUMBER.\@\ \NAME.\@\NOTE},
    notefont=\normalfont\bfseries\boldmath,
    notebraces={}{.},
]{definition}
\declaretheoremstyle[
    spaceabove=\parasep, spacebelow=\parasep,
    postheadspace=.5em,
    headfont=\normalfont\bfseries,
    headpunct={},
    headformat={\NUMBER.\@\NOTE},
    notefont=\normalfont\bfseries\boldmath,
    notebraces={}{.},
]{para}

\renewenvironment{proof}[1][\proofname]{\par
    \pushQED{\qed}%
    \normalfont\trivlist
    \item[\hskip\labelsep\itshape #1\@addpunct{\textbf{.}}]\ignorespaces
}{%
    \popQED\endtrivlist\@endpefalse
}

\makeatother

\declaretheorem[sibling=paragraph, style=para, refname={\S,\S\S}]{para}
\declaretheorem[sibling=paragraph, style=theorem, name=Theorem]{theorem}
\declaretheorem[sibling=paragraph, style=theorem, name=Lemma]{lemma}
\declaretheorem[sibling=paragraph, style=theorem, name=Corollary]{corollary}
\declaretheorem[sibling=paragraph, style=theorem, name=Proposition]{proposition}

\declaretheorem[numbered=no, style=theorem*, name=Theorem]{theorem*}
\declaretheorem[numbered=no, style=theorem*, name=Lemma]{lemma*}
\declaretheorem[sibling=paragraph, style=definition, name=Definition]{definition}
\declaretheorem[sibling=paragraph, style=definition, name=Example]{example}

\declaretheorem[sibling=paragraph, style=definition, name=Remark]{remark}

\numberwithin{equation}{paragraph}

\crefdefaultlabelformat{#2{\upshape#1}#3}
\crefname{figure}{Figure}{Figures}

\crefformat{equation}{#2{\upshape(#1)}#3}
\crefrangeformat{equation}{{\upshape#3(#1)#4--#5(#2)#6}}
\crefmultiformat{equation}{#2{\upshape(#1)}#3}{ and #2{\upshape(#1)}#3}{, #2{\upshape(#1)}#3}{, and~#2{\upshape(#1)}#3}

\crefformat{enumi}{#2{\upshape#1}#3}
\crefrangeformat{enumi}{{\upshape#3#1#4--#5#2#6}}
\crefmultiformat{enumi}{#2{\upshape#1}#3}{ and #2{\upshape#1}#3}{, #2{\upshape#1}#3}{, and~#2{\upshape#1}#3}

\crefformat{section}{#2{\upshape\S#1}#3}
\crefrangeformat{section}{{\upshape#3\S\S#1#4--#5#2#6}}
\crefmultiformat{section}{{\upshape#2\S\S#1#3}}{ and {\upshape#2#1#3}}{, {\upshape#2#1#3}}{, and~{\upshape#2#1#3}}
\crefformat{subsection}{#2{\upshape\S#1}#3}
\crefrangeformat{subsection}{{\upshape#3\S\S#1#4--#5#2#6}}
\crefmultiformat{subsection}{{\upshape#2\S\S#1#3}}{ and {\upshape#2#1#3}}{, {\upshape#2#1#3}}{, and~{\upshape#2#1#3}}
\crefformat{para}{#2{\upshape\S#1}#3}
\crefrangeformat{para}{{\upshape#3\S\S#1#4--#5#2#6}}
\crefmultiformat{para}{{\upshape#2\S\S#1#3}}{ and {\upshape#2#1#3}}{, {\upshape#2#1#3}}{, and~{\upshape#2#1#3}}

\usepackage{enumitem}
\setlist{noitemsep}
\setlist[enumerate]{label=\textnormal{(\roman*)}}

\usepackage{caption}
\DeclareCaptionFont{fira}{\firamedium\firaproportional}
\newcommand{\preparebibliography}{
    \phantomsection
    \addcontentsline{toc}{section}{References}
    \sloppy
    \setstretch{1.1}
    \renewcommand*{\bibfont}{\normalfont\small}
    \setlength{\bibitemsep}{0.2\baselineskip}
}

\usepackage[
    backend=biber,
    style=oxnum,
    giveninits,
    maxbibnames=5,
    maxcitenames=5,
    sorting=nyt
]{biblatex}

\DefineBibliographyStrings{english}{
  withafterword = {with an appendix by},
}

\DeclareFieldFormat[article]{version}{\IfDecimal{#1}{version~}{}#1}
\renewbibmacro*{series+number}{%
    \printfield{series}\isdot%
    \setunit*{\addcomma\addspace}%
    \printfield{number}}

\DeclareFieldFormat{pages}{%
  \iffieldundef{bookpagination}%
    {#1}%
    {\mkpageprefix[bookpagination]{#1}}%
}

        \tikzset{
            symbol/.style={
            draw=none,
            every to/.append style={
            edge node={node [sloped, allow upside down, auto=false]{$#1$}}}
            }
        }

        \usetikzlibrary{decorations.markings}

            \makeatletter
            \tikzcdset{
              open/.code     = {\tikzcdset{hook, circled};},
              closed/.code   = {\tikzcdset{hook, slashed};},
              open'/.code    = {\tikzcdset{hook', circled};},
              closed'/.code  = {\tikzcdset{hook', slashed};},
              circled/.code  = {\tikzcdset{markwith = {\draw (0,0) circle (.375ex);}};},
              slashed/.code  = {\tikzcdset{markwith = {\draw[-] (-.4ex,-.4ex) -- (.4ex,.4ex);}};},
              markwith/.code ={
                \pgfutil@ifundefined%
                {tikz@library@decorations.markings@loaded}%
                {\pgfutil@packageerror{tikz-cd}{You need to say %
                  \string\usetikzlibrary{decorations.markings} to use arrows with markings}{}}{}%
                \pgfkeysalso{/tikz/postaction = {
                  /tikz/decorate,
                  /tikz/decoration={markings, mark = at position 0.5 with {#1}}}
                }
              },
            }
            \makeatother

\usepackage{enotez}
\setenotez{backref=true}

\tikzset{
    labl/.style={anchor=south, rotate=90, inner sep=.5mm}
}
\tikzset{
    lablh/.style={anchor=south, inner sep=.5mm}
}

        \renewcommand{\ss}{{\mathrm{ss}}}
        
        \newcommand{\tr}{{\mathrm{tr}}}
        \newcommand{\pr}{{\mathrm{pr}}}

        \newcommand{\N}{{\ensuremath{{\mathbb N}}}}

        \newcommand{\R}{{\ensuremath{{\mathbb R}}}}

        \newcommand{\Z}{{\ensuremath{{\mathbb Z}}}}
        \newcommand{\Q}{{\ensuremath{{\mathbb Q}}}}
        \newcommand{\G}{{\ensuremath{{\mathbb G}}}}
        \newcommand{\A}{{\ensuremath{{\mathbb A}}}}

        \newcommand{\bA}{{\mathbb A}}
        
        \newcommand{\bC}{{\mathbb C}}

        \newcommand{\bQ}{{\mathbb Q}}

        \newcommand{\bX}{{\mathbb X}}
        
        \newcommand{\bZ}{{\mathbb Z}}

        \newcommand{\cD}{{\mathcal D}}

        \newcommand{\cL}{{\mathcal L}}

        \newcommand{\cO}{{\mathcal O}}

        \newcommand{\abs}[1]{\left|  #1\right| } 

        \makeatletter
        \newcommand{\colim@}[2]{%
          \vtop{\m@th\ialign{##\cr
            \hfil$#1\operator@font colim$\hfil\cr
            \noalign{\nointerlineskip\kern1.5\ex@}#2\cr
            \noalign{\nointerlineskip\kern-\ex@}\cr}}%
        }
        \newcommand{\colim}{%
          \mathop{\mathpalette\colim@{\rightarrowfill@\scriptscriptstyle}}\nmlimits@
        }
        \renewcommand{\varprojlim}{%
          \mathop{\mathpalette\varlim@{\leftarrowfill@\scriptscriptstyle}}\nmlimits@
        }
        \renewcommand{\varinjlim}{%
          \mathop{\mathpalette\varlim@{\rightarrowfill@\scriptscriptstyle}}\nmlimits@
        }
        \makeatother

        \newcommand{\git}{
          \mathchoice{/\mkern-6mu/}
            {/\mkern-6mu/}
            {/\mkern-5mu/}
            {/\mkern-5mu/}}

        \DeclareMathOperator{\im}{im}

        \newcommand{\ev}{ev}
        
        \newcommand{\Hom}{\mathrm{Hom}}

        \newcommand{\Aut}{\mathrm{Aut}}
        
        \newcommand{\GL}{\mathrm{GL}}
        
        \newcommand{\Ad}{\mathrm{Ad}}

        \DeclareMathOperator{\Spec}{\mathrm{Spec}}

        \newcommand{\st}{\vert\ }

        \newcommand{\Filt}{\mathrm{Filt}}
        \newcommand{\Grad}{\mathrm{Grad}}
        \let\ev\relax
        \newcommand{\ev}{\mathrm{ev}}
        \newcommand{\gr}{\mathrm{gr}}

        \newcommand{\norm}[1]{\lVert #1\rVert}

        \DeclareMathOperator{\Bun}{Bun}

        \DeclareMathOperator{\wt}{\mathrm{wt}}

        \DeclareMathOperator{\qfilt}{\Q\,\mathrm{-}\Filt}
        
        \renewcommand{\phi}{\varphi}

        \newcommand{\CL}{\mathrm{CL}}

        \makeatletter 
        \DeclareFontEncoding{LS1}{}{}
        \DeclareFontSubstitution{LS1}{stix}{m}{n}
        \DeclareMathAlphabet{\mathcaldos}{LS1}{stixscr}{m}{n}
        \makeatother

\let\Cref\cref
        
\title{Stability and disconnected groups}
\author{Andres Fernandez Herrero\and Andrés Ibáñez Núñez}
\date{}

\begin{document}
\initlengths

\maketitle
\vspace{20pt}
\begin{abstract}
We study the notion of semistability for principal bundles over curves with possibly disconnected reductive structure group. We establish a new characterization of the behavior of semistability under change of group, novel even in the connected case, and prove that all existing notions of semistability are equivalent, thus settling a question by Biswas-Gomez. The key ingredients for our results include a study of cocharacters and characters of disconnected linear algebraic groups, and an extension of the recursive description of Kirwan stratifications in Geometric Invariant Theory to the case of disconnected groups.
\end{abstract}
\vspace{20pt}
\setcounter{secnumdepth}{3}
\setcounter{tocdepth}{2}
{
    \hypersetup{linkcolor=black}
    \setstretch{1.1}
    \tableofcontents
}

\newpage

\section{Introduction}

\addtocounter{subsection}{1}

In this paper we provide a treatment of semistability of $G$-bundles on a smooth projective connected curve $C$ over an algebraically closed field $k$, where $G$ is a possibly disconnected reductive algebraic group. We define a notion of semistability in terms of Hilbert-Mumford weights of the determinant line bundle on the corresponding stack of $G$-bundles $\Bun_G$ (\Cref{defn: semistability bundles}), which can be equivalently described in a more classical way using parabolic reductions (see \ref{para: semistability in terms of parabolic reductions}). One of our main results is a description of how semistability behaves under change of group.

\begin{theorem}[{($=$ \Cref{thm: semistability under change of group})}] \label{thm: intro thm semistability under change of group}
    Suppose that the algebraically closed ground field $k$ has characteristic $0$. Let $f: G \to H$ be a homomorphism of linearly reductive affine algebraic groups over $k$. Let $E$ be a $G$-bundle on $C$. Then, the following hold:
    \begin{enumerate}
        \item If the rational degree of $E$ (see \ref{para: rational degree of G-bundle}) is not $f$-adapted (as in \Cref{definition: adapted degree}), then the associated $H$-bundle $E(H)$ is not semistable.
        \item If the rational degree of $E$ is $f$-adapted and $E$ is semistable, then the associated $H$-bundle $E(H)$ is semistable.
    \end{enumerate}
\end{theorem}

We encourage the reader to look at \Cref{example: rational degrees for direct sum homomorphism} and \Cref{example: semistability under direct sum homomorphism} to get a sense of \Cref{thm: intro thm semistability under change of group} in a simple toy case.

To our knowledge, the general result in \Cref{thm: intro thm semistability under change of group} is not found in the literature even in the case when the groups $G$ and $H$ are connected (see \Cref{example: the case of connected groups}, where we spell out \Cref{thm: intro thm semistability under change of group} in this special case).
Some precursors of \Cref{thm: intro thm semistability under change of group} in the case when both $H$ and $G$ are connected include:
\begin{itemize}
    \item the case when the image of the center $Z(G)$ is contained in $Z(H)$ is proven in Ramanathan's thesis \cite[Prop. 3.17]{ramanathan_moduli_principal_bundles} and the article by Ramanan and Ramanathan \cite[Thm. 3.18]{ramanan-ramanathan}; and
    \item the case where $f:G \to H$ is the inclusion of a Levi subgroup of a parabolic $P \subset H$ is proven by Frățilă in \cite[Lemmas 2.14 and 2.15]{fratila_g_bundles_elliptic_curves}.
\end{itemize}

There are several notions of semistability of $G$-bundles for disconnected groups in the literature. These include:
\begin{itemize}
    \item Hilbert-Mumford semistability, the notion of semistability we use in this paper. 
    \item ad-semistability, as in \cite{_Atiyah_TheYangMillsEquationsoverRiemannSurfaces, moduli-disconnected-groups}.
    \item Einstein-Hermite semistability, as in \cite{biswas-gomez-kahler-einstein}.
\end{itemize}
We refer the reader to \Cref{section: comparison notions semistability} for more details on the definitions and a more complete list of references and precursors. Using \Cref{thm: intro thm semistability under change of group}, we show that all existing notions of semistability for $G$-bundles are equivalent (see \Cref{prop: survey all semistability conditions}), thus settling a question by \textcite{biswas-gomez-kahler-einstein}.

\begin{remark}
    The arguments in this paper can be modified using the approach in \cite{herrero2023meromorphic} to prove that \Cref{thm: intro thm semistability under change of group}(ii) holds in the context of $G$-Higgs bundles, $G$-bundles with connections, and more generally $G$-bundles with a $t$-connection. This can be used as in \cite{herrero2023meromorphic} to develop the theory of moduli spaces of $G$-Higgs bundles and $G$-bundles with connection in the setting of disconnected groups. We have decided not to include such variations this in the paper in the interest of space, and leave the details to the interested reader.
\end{remark}

To prove \Cref{thm: intro thm semistability under change of group}, we need to use certain results on (co)characters of disconnected groups (\cref{theorem: main intro thm}) and Kirwan stratifications for disconnected algebraic groups (\cref{thm: kirwan stratification intro theorem}). Even though these results should not be surprising to experts, they seem to be absent in the literature, and so we have decided to include a complete treatment with proofs.

The following \Cref{theorem: main intro thm} collects some facts on (co)characters of disconnected groups, which we prove in Section \ref{section: weyl groups and cocharacters}. These results have been applied in the appendix of \cite{dervan2024stabilityconditionsgeometricinvariant}, which cites a previous version of this paper.
\begin{theorem}\label{theorem: main intro thm} 
Let $G$ be a possibly disconnected algebraic group over an algebraically closed field $k$, with maximal torus $T$.
\begin{enumerate}
\item The Weyl group $W=N_G(T)/Z_G(T)$ of $G$ is isomorphic to the Weyl group of its reductive quotient $G/R_\mathrm{u}(G)$.
\item The restriction map $\mathbb{X}^*(G)_\bQ\to \mathbb{X}^*(T)_\bQ$ on rational characters induces an isomorphism $\mathbb{X}^*(G)_\bQ\cong \mathbb{X}^*(T)^W_\bQ$ to the Weyl-invariant part of $\mathbb{X}^*(T)_\bQ$.

\item If $G$ is reductive with centre $Z(G) \subset G$, then the restriction map $\mathbb{X}^*(G)_{\mathbb{Q}} \to \mathbb{X}^*(Z(G))_\bQ$ on rational characters and the inclusion map $\mathbb{X}_*(Z(G)) \to \mathbb{X}_*(T)^W$ on integral cocharacters are isomorphisms.
\end{enumerate}
\end{theorem}

We note that \Cref{theorem: main intro thm} is well-known in the case of connected algebraic groups, see \cite[\S 24.1, Prop. B]{humphreys_linear_algebraic_groups} and \cite[\S 26.2, Cor. B]{humphreys_linear_algebraic_groups}.
However, to our knowledge, there is no explicit treatment of \Cref{theorem: main intro thm} in the literature; one of the purposes of this paper is to provide a citable reference for these facts.

In Section \ref{section: kirwan stratification}, we use \Cref{theorem: main intro thm} to extend the classical description of the Kirwan stratification to the case of disconnected  reductive groups. We assume that the algebraically closed ground field $k$ has characteristic $0$, and that the linear algebraic group $G$ is reductive. While Kirwan's original construction deals with an action of $G$ on a smooth projective variety, we work in a more general relative setting. We fix a proper $G$-equivariant morphism $f: X \to S$ between finite-type separated schemes over $k$ endowed with actions of $G$. We assume that the quotient stack $S/G$ has a good moduli space $\pi: S/G \to S\git G$ as in  \textcite{_Alper_GoodmodulispacesforArtinstacks} (this is automatically satisfied if $S$ is affine). To define the instability stratification, we also fix the following data:
\begin{itemize}
    \item A $G$-equivariant rational line bundle $\cL$ on $X$ that is $f$-ample.
    \item A $W$-invariant rational inner product $(-,-)$ on the vector space $\mathbb{X}_*(T)_\bQ$.
\end{itemize}
We refer the reader to Section \ref{section: kirwan stratification} for more details on the following result, which we state in an informal way for the purposes of this introduction.

\begin{theorem}[($=$ {\Cref{thm: description of centres})}] \label{thm: kirwan stratification intro theorem}
With assumptions as above, there is a $G$-equivariant instability stratification of $X$ by locally closed subschemes (\Cref{thm: kirwan stratification}). Each stratum can be recursively described in terms of the semistable locus $X^{\mathrm{ss}}_{\alpha}$ of a locus of fixed points $X_{\alpha} \subset X$ with respect to an explicitly described linearization for a Levi subgroup of $G$ (\Cref{thm: description of centres}).
\end{theorem}

We remark that, in the case when $S = \Spec(k)$, the statement of \Cref{thm: description of centres} is claimed in \cite[Remark 2.21]{_Kirwan_Cohomologyofquotientsinsymplecticandalgebraicgeometry}. However, a proof was not provided, and it seems that one would need a consequence of \Cref{theorem: main intro thm} to proceed.

\cref{thm: kirwan stratification intro theorem} is a crucial ingredient in the proof of \Cref{thm: intro thm semistability under change of group}, where it is applied to the action of $G$ on a certain flag variety $H/Q$. \cref{thm: kirwan stratification intro theorem} is also of independent interest, and we establish it in a more general relative setting with the expectation that it will be useful in other applications.

\medskip
 \noindent \textbf{Acknowledgements.}
 We would like to thank Brian Conrad and Dragoș Frățilă for helpful email exchanges. We are grateful to Tom\'as G\'omez, Tasuki Kinjo, Frances Kirwan and Alfonso Zamora for useful comments. 

 \medskip
 \noindent \textbf{Notations and conventions.}
 We write $\N=\Z_{\geq 0}$ for the set of natural numbers. For an algebraic group $G$ over a field $k$ acting on a $k$-scheme $X$, we write $X/G$ for the quotient stack, without the usual bracket decorations. Our reductive groups are not assumed to be connected.

Given an algebraic group $H$ over a field $k$, we denote by $\mathbb{X}^*(H):= \Hom(H, \mathbb{G}_{m,k})$ the group of characters and $\mathbb{X}_*(H)$ the set of cocharacters. The $\Q$-vector space of rational characters is denoted $\mathbb{X}^*(H)_\bQ := \mathbb{X}^*(H) \otimes_{\mathbb{Z}} \mathbb{Q}$. Even though $\mathbb{X}_*(H)$ may not be a group, we can still define the set of rational cocharacters as 
\[\mathbb{X}_*(H)_\bQ=\{\lambda/n  \,\vert\,  \lambda\in \mathbb{X}_*(G),\ n\in \Z_{>0}\}/\sim,\]
where $\lambda_1/n_1\sim \lambda_2/n_2$ if $\lambda_1^{n_2}=\lambda_2^{n_1}$.

If $H$ is a torus, then $\mathbb{X}_*(H)$ is naturally equipped with the structure of a group, and we have $\mathbb{X}_*(H)_\bQ = \mathbb{X}_*(H) \otimes_{\mathbb{Z}} \mathbb{Q}$. In this case we have a canonical identification $\mathbb{X}_*(H)_\bQ = \mathbb{X}^*(H)^{\vee}_\bQ$, and we denote by $\langle -, - \rangle: \mathbb{X}_*(H)_\bQ \otimes_{\mathbb{Q}} \mathbb{X}^*(H)_\bQ \to \mathbb{Q}$ the induced nondegenerate pairing.

Given a vector space $V$ equipped with an action of a finite group $W$, we denote by $V^W$ and $V_W$ the vector spaces of invariants and coinvariants, respectively.

\section{Weyl groups, cocharacters and rational characters} \label{section: weyl groups and cocharacters}

In this section we prove \Cref{theorem: main intro thm} from the introduction. We fix an algebraically closed field $k$ and a smooth affine algebraic group $G$ over $k$. For standard facts about linear algebraic groups and reductive groups we refer the reader to \cite{_Milne_Algebraicgroupsthetheoryofgroupschemesoffinitetypeoverafield,_Conrad_PseudoreductiveGroups}.

\subsection{Weyl groups}

\begin{para}
   We fix a maximal torus $T$ of $G$. Let us denote by $N_G(T)$ the normalizer of $T$ in $G$ and by $Z_G(T)$ the centralizer of $T$ in $G$. The \emph{Weyl group} (of $G$ with respect to $T$) is the quotient $W=W(G,T)=N_G(T)/Z_G(T)$. The Weyl group embeds in the automorphism group scheme $\underline\Aut(T)$ of the torus $T$, and it is thus a finite constant group. 
\end{para}

\begin{para}
    Let $U$ be the unipotent radical of $G$ and let $R=G/U$ denote its \emph{Levi quotient}. The group $G$ is said to be \emph{reductive} if $U=1$. We denote by $T'$ the image of $T$ in $R$, which is a maximal torus of $R$ isomorphic to $T$, and let $W'=W(R, T')$ be the associated Weyl group.
\end{para}

 The unipotent radical does not play a role in the formation of the Weyl group:

\begin{proposition}\label{proposition: Weyl group only sees reductive}
    There is a canonical isomorphism $W\cong W'$.
\end{proposition}

\begin{proof}
    We have a diagram
    \[\begin{tikzcd}
        & 1 & 1 & 1 \\
        1 & {Z_U(T)} & {N_U(T)} & {N_U(T)/Z_U(T)} & 1 \\
        1 & {Z_G(T)} & {N_G(T)} & W & 1 \\
        1 & {Z_R(T')} & {N_R(T')} & {W'} & 1  
        \arrow[from=2-1, to=2-2]
        \arrow[from=2-2, to=2-3]
        \arrow[from=2-3, to=2-4]
        \arrow[from=2-4, to=2-5]
        \arrow[from=3-1, to=3-2]
        \arrow[from=3-2, to=3-3]
        \arrow[from=3-3, to=3-4]
        \arrow[from=3-4, to=3-5]
        \arrow[from=4-1, to=4-2]
        \arrow[from=4-2, to=4-3]
        \arrow[from=4-3, to=4-4]
        \arrow[from=4-4, to=4-5]
        \arrow[from=1-2, to=2-2]
        \arrow[from=2-2, to=3-2]
        \arrow[from=3-2, to=4-2]
        \arrow[from=1-3, to=2-3]
        \arrow[from=2-3, to=3-3]
        \arrow[from=3-3, to=4-3]
        \arrow[from=1-4, to=2-4]
        \arrow[from=2-4, to=3-4]
        \arrow[from=3-4, to=4-4]
    \end{tikzcd}\]
    where all rows and columns are exact. To conclude, it is enough to show that:
    
    \begin{enumerate}[label=(\arabic*)]
        \item \label{item: first proof of prop 1} $N_G(T)\to N_R(T')$ is surjective.
        \item \label{item: second proof of prop 1} The inclusion $Z_U(T) \hookrightarrow N_U(T) $ is an isomorphism.
    \end{enumerate}
    
    For \ref{item: first proof of prop 1} we first note that the preimage $V$ of $T'$ along $G\to R$ is the semidirect product $V=U\rtimes T$. Indeed, $V=U\cdot T$, $T$ normalises $U$ and, since $T$ is a torus and $U$ is unipotent, we have $T\cap U=1$. Let $r\in N_R(T')(k)$ and let $g\in G(k)$ be a lift of $r$ to $G$. Then $gTg^{-1}$ is a maximal torus of $V$, so there exists $u\in U(k)$ such that $gTg^{-1}=uTu^{-1}$, by conjugacy of maximal tori. We have that $u^{-1}g$ is in $N_G(T)$ and that it is a preimage of $r$ along $N_G(T)\to N_{R}(T')$. This shows \ref{item: first proof of prop 1}.
    
    For \ref{item: second proof of prop 1}, let $A$ be a commutative $k$-algebra and let $g\in N_U(T)(A)$. For every $A$-algebra $B$ and every $t\in T(B)$, we have that $tg^{-1}_{\vert {B}}t^{-1}\in U(B)$, because $U$ is normal, and $g_{\vert B} t g^{-1}_{\vert B}\in T(B)$, because $g$ normalises $T$. Therefore $g_{\vert B} t g^{-1}_{\vert B} t^{-1}\in T(B)\cap U(B)$. Since $T\cap U$ is trivial, we must have $g\in Z_U(T)(A)$. 
\end{proof}

\subsection{Cocharacters}
In the reductive case, central cocharacters of $G$ coincide with the Weyl-invariant cocharacters of the maximal torus $T$:

\begin{proposition}\label{theorem: cocharacters and Weyl group}
    Suppose that $G$ is reductive with centre $Z(G) \subset G$. Then the natural map $\mathbb{X}_*(Z(G))\to \mathbb{X}_*(T)^W$ is an isomorphism.
\end{proposition}

We first give a self-contained argument using dynamic techniques. Then, below, we explain how the result can also be deduced from the connected case.

\begin{proof}
The reduced identity component $Z(G)^\circ_{\mathrm{red}}$ of the centre $Z(G)$ of $G$ is a torus contained in $T$, so we have an injection
\[\mathbb{X}_*(Z(G))=\mathbb{X}_*(Z(G)^\circ_{\mathrm{red}})\hookrightarrow \mathbb{X}_*(T)^W.\]

Let $\lambda\in \mathbb{X}_*(T)^W$. We want to prove that $\lambda\colon \G_{m,k}\to T$ factors through $Z(G)$ or, equivalently, that the centraliser $L_G(\lambda)$ of $\lambda$ in $G$ equals the whole of $G$.

We denote $P_G(\lambda)$ and $U_G(\lambda)$ the subgroups of $G$ defined functorially by

\[P_G(\lambda)(A)=\{g\in G(A)\st \lim_{t\to 0}\lambda(t)g\lambda(t)^{-1} \text{ exists in } G\},\]
\[U_G(\lambda)(A)=\{g\in P_G(\lambda)(A)\st \lim_{t\to 0}\lambda(t)g\lambda(t)^{-1}=1\},\]
for any $k$-algebra $A$. By existence of the limit we mean that the corresponding morphism $\G_{m,A}\to G$ extends to a morphism $\A^1_A\to G$. The subgroups $P_G(\lambda)$ and $U_G(\lambda)$ are smooth because $G$ is \cite[Proposition 2.1.8]{_Conrad_PseudoreductiveGroups}. Furthermore, $U_G(\lambda)$ is connected (even if $G$ is not) because it retracts to a point, and it is also unipotent \cite[Lemma~2.1.5]{_Conrad_PseudoreductiveGroups}.

The result follows from the following:

\medskip

\noindent \textit{Claim.} \textit{We have $P_G(\lambda) = G$.}

\medskip
Indeed, assume the claim. Note that $P_G(\lambda)=U_G(\lambda)\rtimes L_G(\lambda)$, where $U_G(\lambda)$ is smooth, connected and unipotent. Since $G$ is reductive, $U_G(\lambda)=1$, and thus $G=L_G(\lambda)$, as desired.

We are left to prove the claim. Since $P_G(\lambda)$ and $G$ are smooth, it suffices to show that for any $g \in G(k)$ we have $g \in P_G(\lambda)(k)$. Denote $\lambda^g=g\lambda g^{-1}$. We first show that there are $p, h \in P_G(\lambda)$ such that $\lambda^{pgh} \in \mathbb{X}_*(T)$. Since the intersection of the two parabolic subgroups $P_{G^\circ}(\lambda), P_{G^\circ}(\lambda^g) \subset G^\circ$ contains a maximal torus of $G^\circ$, there exist $h_1\in P_{G^\circ}(\lambda^g)(k)$ and $h_2\in P_{G^\circ}(\lambda)(k)$ such that $\lambda^{h_1g}$ and $\lambda^{h_2}$ commute. Since $P_{G^\circ}(\lambda^g)=gP_{G^\circ}(\lambda)g^{-1}$, there are $b,h\in P_{G^\circ}(\lambda)(k)$ such that $\lambda^{bgh}$ and $\lambda$ commute (just take $b=h_2^{-1}$ and $h=g^{-1}h_1g$). Let $T_1$ be a maximal torus containing both $\lambda^{bgh}$ and $\lambda$. Both $T_1$ and $T$ are maximal tori of $L_G(\lambda)$, so there is $u\in L_G(\lambda)(k)$ such that $uT_1 u^{-1}=T$. Let $p= ub$, then $\lambda^{pgh}\in \mathbb{X}_*(T)$.

Now, set $v=pgh$. We have that $v^{-1}Tv$ is a maximal torus containing $\lambda$, so there is $r\in L_G(\lambda)(k)$ such that $rTr^{-1}=v^{-1}Tv$. Then $vr\in N_G(T)(k)$. Since $\lambda$ is fixed by $W$, we have $\lambda=\lambda^{vr}=(\lambda^r)^v=\lambda^v$, and hence $v=pgh\in L_G(\lambda)(k)$. Since $p,h\in P_G(\lambda)(k)$, this implies that $g\in P_G(\lambda)(k)$ as well.
\end{proof}

\begin{para}[An alternative argument]
If one is willing to accept \Cref{theorem: cocharacters and Weyl group} in the connected case, one can also argue as follows in the general case.
Since $G^\circ$ is connected reductive, the automorphism group can be written as $\Aut(G^\circ)=G^\circ/Z(G^\circ) \rtimes \operatorname{Out}(G^\circ)$, where $\operatorname{Out}(G^\circ)$ preserves pinning (in particular, the torus $T$) \cite[Prop. 7.1.6]{conrad_reductive_group_schemes}. 
Thus, if $g\in G(k)$, then conjugation by $g$ in $G^\circ$ is the same as conjugation by $ah$, with $a\in G^\circ(k)$ and $h\in N_G(T)$. If $\lambda\in \mathbb{X}_*(T)^W$, then $\lambda\in \mathbb{X}_*(T)^{W(G^\circ,T)}=\mathbb{X}_*(Z(G^\circ))$ by the connected case. Then $\lambda^g=(\lambda^h)^a=\lambda^a=\lambda$. Since this is true for all $g\in G(k)$, we have $\lambda\in \mathbb{X}_*(Z(G))$, as desired. 
\end{para}

\subsection{Rational characters}

\begin{para} 
   Every rational character $\chi$ of $G$ restricts to a rational character $\chi\vert_T$ of $T$, and moreover, this restriction is $W$-invariant. In this way we get a map $\mathbb{X}^*(G)_\bQ\to \mathbb{X}^*(T)_\bQ^W$.
\end{para}

 Rational characters of $G$ can be described in terms of those of $T$ and the Weyl group:
 
\begin{proposition}\label{theorem: rational characters and Weyl group}
    The natural map $\mathbb{X}^*(G)_\bQ\to \mathbb{X}^*(T)_\bQ^W$
    is an isomorphism. Moreover, if $G$ is reductive, then $\mathbb{X}^*(G)_\bQ\cong \mathbb{X}^*(Z(G))_\bQ$,
    where $Z(G)$ is the centre of $G$.
\end{proposition}

\begin{proof}
First we note that it is sufficient to show the result in the case when $G$ is reductive. Indeed, if the result was true in the reductive case, then for an arbitrary $G$ we would have
\[\mathbb{X}^*(G)_\bQ=\mathbb{X}^*(R)_\bQ\cong\mathbb{X}^*(T')_\bQ^{W'}=\mathbb{X}^*(T)_\bQ^W,\]
because every character $G\to \G_{m,k}$ factors through $R$, because the result is true for $R$, and by \Cref{proposition: Weyl group only sees reductive}. Therefore, we may assume without loss of generality that $G$ is reductive.

Recall that if a finite group $S$ acts on a finite dimensional $\Q$-vector space $V$, then there is a canonical bijection $(V^S)^\vee=(V^\vee)^S$. Therefore, dualizing the equality $\mathbb{X}_*(Z(G))_\bQ=\mathbb{X}_*(T)_\bQ^W$ from \Cref{theorem: cocharacters and Weyl group}, we get $\mathbb{X}^*(Z(G))_\bQ = \mathbb{X}^*(T)_\bQ^{W}$. In view of this, for the identification $\mathbb{X}^*(G)_\bQ = \mathbb{X}^*(T)_\bQ^{W}$, it suffices to show that the restriction morphism $\varphi: \mathbb{X}^*(G)_\bQ \to \mathbb{X}^*(Z(G))_\bQ$ is bijective.

To see injectivity of $\varphi$, we use the following commutative diagram
\[\begin{tikzcd}
	{\mathbb{X}^*(G)_\bQ} & {\mathbb{X}^*(G^\circ)_\bQ} \\
	{\mathbb{X}^*(Z(G))_\bQ} & {\mathbb{X}^*(Z(G^\circ))_\bQ}
	\arrow[from=1-1, to=1-2]
	\arrow[from=1-1, to=2-1]
	\arrow[from=1-2, to=2-2]
	\arrow[from=2-1, to=2-2]
\end{tikzcd}\]

In view of this diagram, it suffices to show that both $\mathbb{X}^*(G^\circ)_\bQ \to \mathbb{X}^*(Z(G^\circ))_\bQ$ and $\mathbb{X}^*(G)_\bQ \to \mathbb{X}^*(G^\circ)_\bQ$ are injective. 

Since $G^\circ$ is connected and reductive, we have the central isogeny $Z(G^\circ)\to G^\circ/\cD(G^\circ)$, where $\cD(G^\circ)$ is the derived subgroup. This shows that $\mathbb{X}^*(G^\circ)_\bQ \to \mathbb{X}^*(Z(G^\circ))_\bQ$ is injective. 

To see $\mathbb{X}^*(G^\circ)_\bQ \to \mathbb{X}^*(Z(G^\circ))_\bQ$ is injective. Let $\chi\in \mathbb{X}^*(G)_\bQ$ such that $\chi|_{G^\circ}=0$. After scaling, we may assume that $\chi$ is represented by a character $\chi\colon G\to \G_{m,k}$ such that $\chi\vert_{G^\circ}$ is trivial. Note that $G/\G^{\circ}$ is a finite group; let $N$ denote its order. For all $g\in G(k)$, $\chi^{N}(g)=\chi(g^{N})=1$, because $g^{N}$ is in $G^\circ$. In additive notation, $N\chi=0$, which means that $\chi =0$. This shows injectivity of $\mathbb{X}^*(G^\circ)_\bQ \to \mathbb{X}^*(Z(G^\circ))_\bQ$, and therefore concludes the proof of injectivity of $\varphi$.

For surjectivity of $\varphi$, we choose an embedding $G\hookrightarrow \GL_{n,k}$. Let $T'$ be the maximal torus of $Z(G)$. We have $G \hookrightarrow Z_{\GL_{n,k}}(T')$. Now, there are positive integers $n_1,\ldots, n_l$ such that $Z_{\GL_{n,k}}\cong\GL_{n_1,k}\times \cdots \GL_{n_l,k}$. Let $\G_{m,k}^l\subset Z_{\GL_{n,k}}$ be the central torus. We have the diagram 
\[\begin{tikzcd}
	{\mathbb{X}^*(G)_\bQ} && {\mathbb{X}^*(\GL_{n_1,k}\times\cdots\times\GL_{n_l,k})_\bQ} \\
	{\mathbb{X}^*(T')_\bQ} && {\mathbb{X}^*(\G_{m,k}^l)_\bQ}
	\arrow["\phi"{swap},from=1-1, to=2-1]
	\arrow[from=1-3, to=1-1]
	\arrow["\sim"{marking, allow upside down, yshift=0.4em}, no head, from=1-3, to=2-3]
	\arrow[two heads, from=2-3, to=2-1]
\end{tikzcd}\]
where the right vertical arrow is an isomorphism and the lower horizontal arrow is surjective. This implies surjectivity of $\phi$.

The last statement in the theorem follows because
\[\mathbb{X}^*(Z(G))_\bQ=\mathbb{X}_*(Z(G))_\bQ^\vee=(\mathbb{X}_*(T)_\bQ^W)^\vee=\mathbb{X}^*(T)_\bQ^W\]
by \Cref{theorem: cocharacters and Weyl group}.
\end{proof}

\section{Kirwan stratifications for disconnected groups} \label{section: kirwan stratification}

\addtocounter{subsection}{1}

    In this section, we work over an algebraically closed field $k$ of characteristic $0$. We start by recalling Kirwan $\Theta$-stratification \cite{_Kirwan_Cohomologyofquotientsinsymplecticandalgebraicgeometry} for certain quotient stacks $X/G$. While Kirwan's original construction assumes $X$ to be a smooth projective variety, we explain how to relax this hypothesis using the theory of $\Theta$-stratifications as developed in \cite{_HalpernLeistner_Onthestructureofinstabilityinmodulitheory}. Then, using \Cref{theorem: rational characters and Weyl group}, we will give an alternative description of the stratification that is crucial in applications. This description used to be available only in the case when $G$ is connected, and it makes crucial use of \Cref{theorem: main intro thm}.

    \begin{para}[Setup and notation]
    We fix a reductive group $G$ over $k$ and a finite-type separated scheme $X$ over $k$ endowed with an action of $G$. We fix a maximal torus $T \subset G$ with Weyl group $W$. Since there will be no ambiguity, in this section we denote $L_\lambda=L_G(\lambda)$ and $P_\lambda=P_G(\lambda)$ for $\lambda\in \mathbb{X}_*(G)_\bQ$ a rational one-parameter subgroup. 
\end{para}

\begin{para}[Stability data]
    Kirwan's stratification will depend on two additional data, that we fix for the rest of this section, namely
    \begin{enumerate}
        \item a rational line bundle $\cL$ on the quotient stack $X/G$ or, equivalently, a rational line bundle $\cL\vert_{X}$ on $X$ with a $G$-equivariant structure, and
        \item a norm on cocharacters of $G$, that is, a $W$-invariant rational inner product $(-,-)$ on the space $\mathbb{X}_*(T)_\bQ$ of cocharacters of $T$.
    \end{enumerate}
\end{para}

\begin{para}[Assumptions]
    We assume that there is another finite-type separated scheme $S$ over $k$ with a $G$-action, and a proper $G$-equivariant morphism $f\colon X\to S$ such that $\cL_{\vert X}$ is $f$-ample. Further, we assume that $S/G$ has a good moduli space $\pi: S/G \to S\git G$ as in  \textcite{_Alper_GoodmodulispacesforArtinstacks}, and that $S\git G$ is a scheme. This condition is automatic if $S=\Spec A$ is affine, in which case the good moduli space is $S\git G=\Spec A^G$.
\end{para}

\subsection{Semistable locus and the Hilbert--Mumford criterion}

\begin{para}[Rational filtrations]
    Suppose that $x\in X(k)$ is a point. For an integral cocharacter $\lambda\in \mathbb{X}_*(G)$, we say that \emph{the limit $\lim_{t\to 0} \lambda(t)x$ exists in $X$} if the morphism
    \[\G_{m,k}\to X\colon \; \; t \mapsto \lambda(t)x\]
    extends to a morphism $\A^1_k\to X$. Since $X$ is separated, the extension is unique if it exist. For a rational cocharacter $\lambda\in \mathbb{X}_*(G)_\bQ$, we say that \emph{the limit $\lim_{t\to 0} \lambda(t)x$ exists in $X$} if for some (and thus for any) positive integer $r$ such that $\lambda^r$ is integral we have that the limit $\lim_{t\to 0} \lambda^r(t)x$ exists in $X$. We define the set of \emph{rational filtrations of $x$ in $X/G$} to be
    \[\qfilt(X/G,x)=\left\{\left.\lambda \in \mathbb{X}_*(G)_\bQ\right\vert \lim_{t\to 0} \lambda(t)x \text{ exists in } X\right\}/\sim,\]
    where $\sim$ is the equivalence relation that identifies $\lambda$ and $\lambda^g$ if $g\in P_\lambda(k)$. 
\end{para}

\begin{para}[Hilbert--Mumford weight]
Now suppose that $\lambda\in \mathbb{X}_*(G)$ is a cocharacter such that $\lim_{t\to 0} \lambda(t)x=y$ exists in $X$. The point $y$ gives a morphism $y\colon \Spec k\to X$ that is equivariant with respect to the homomorphism $\lambda\colon \G_{m,k}\to G$, so it induces a map $u\colon \mathrm{B}\G_{m,k}=\Spec k/\G_{m,k}\to X/G$ between quotient stacks. The pullback $\cL_{y,\lambda}=u^*\cL$ is a line bundle on $\mathrm{B}\G_{m,k}$, and thus a one-dimensional representation of $\G_{m,k}$. The \emph{Hilbert--Mumford weight} $\mathrm{m}(x,\lambda)$ is defined to be the opposite of the weight of this representation:
\[\mathrm{m}(x,\lambda)=-\mathrm{wt} \left(\cL_{y,\lambda}\right).\]

For a positive integer $n$, we have $\mathrm{m}(x,\lambda^n)=n\mathrm{m}(x,\lambda)$. Also, if $g\in P_\lambda(k)$, then $\mathrm{m}(x,\lambda^g)=\mathrm{m}(x,\lambda)$. Thus $\mathrm{m}(x,-)$ extends to a map $\mathrm{m}(x,-)\colon \qfilt(X/G,x)\to \Q$.
\end{para}

\begin{definition}[Semistable locus] \label{defn: semistable locus}
    We say that $x \in X(k)$ is \emph{semistable} (with respect to $\cL$ and relative to $S\git G$) if there is an affine open subscheme $U \subset S\git G$ with corresponding preimage $X_U \subset X$ in $X$, a positive integer $k>0$, and a $G$-invariant section $s \in \Gamma(X_U, \cL^{\otimes k}\vert_X)^G$ such that the open complement $(X_U)_s \subset X_U$ of its vanishing locus is affine and $x \in (X_U)_s(k)$.

    A semistable point $x \in X(k)$ is said to be \emph{polystable} if the orbit morphism $G \to X$ given by $g \mapsto g \cdot x$ has closed image. A polystable point $x$ is called \emph{stable} if the $G$-stabilizer $G_x$ of $x$ is finite.
\end{definition}

\begin{para}[Good moduli space for the semistable locus]\label{para-good-moduli-ss-locus}
Under our assumptions, the set of semistable points as in \Cref{defn: semistable locus} are the $k$-points of a (unique) $G$-equivariant open subscheme $X^{\mathrm{ss}} \subset X$, and the corresponding quotient stack $X^{\mathrm{ss}}/G$ admits a projective good moduli space over $S\git G$ given by $\mathrm{Proj}_{S\git G}\left(\left(\bigoplus_{n \in \mathbb{N}} (\pi \circ f)_*\cL|_X^{\otimes n}\right)^G\right)$, see \cite[Theorem 5.6.1]{_HalpernLeistner_Onthestructureofinstabilityinmodulitheory}.
\end{para}

\begin{theorem}[Hilbert--Mumford criterion] \label{thm: hilbert-mumford criterion}
    The following statements hold.
    \begin{enumerate}
        \item A point $x\in X(k)$ is semistable if and only if, for all $\lambda\in \qfilt(X/G,x)$, we have $\mathrm{m}(x,\lambda)\leq 0$.
        \item A point $x \in X(k)$ is polystable if and only if we have $\mathrm{m}(x,\lambda)\leq 0$ for all $\lambda$ as in part (i), with equality if and only if there is some $g \in P_{\lambda}(k)$ such that $g\im(\lambda)g^{-1} \subset G_x$.
        \item A point $x \in X(k)$ is stable if and only if we have $\mathrm{m}(x,\lambda)\leq 0$ for all $\lambda$ as in part (i), with equality if and only if $\lambda =0$.
    \end{enumerate}
\end{theorem}
\begin{proof}
This is the main result in \cite{Gilbrandsen-Halle-Hulek-relative-hm}.
\end{proof}

\subsection{Fixed points and attractors}

\begin{para}[Fixed points]
    For a one-parameter subgroup $\lambda\colon \G_{m,k}\to G$, we have an induced $\G_{m,k}$-action on $X$. We denote $X^{\lambda,0}$ the fixed point locus for that action. If $r\in \Z_{>0}$, then $X^{r\lambda,0}=X^{\lambda,0}$, and thus $X^{\lambda,0}$ makes sense if $\lambda$ is a rational one-parameter subgroup.
\end{para}

\begin{para}[Attractors]
The \emph{attractor} $X^{\lambda,+}$ is the algebraic space whose $R$-points for a $k$-scheme $R$ are given by
\[X^{\lambda,+}(R)=\Hom^{\G_{m,k}}(R\times \A^1,X),\]
that is, a map $R\to X^{\lambda,+}$ is a $\G_{m,k}$-equivariant map $R\times \A^1\to X$, where the action on $R\times \A^1$ is by scaling and the action on $X$ is the one induced by $\lambda$. The attractor $X^{\lambda,+}$ is represented by an algebraic space, see \cite[Prop. 1.4.1]{_HalpernLeistner_Onthestructureofinstabilityinmodulitheory} and \cite{_Drinfeld_OnalgebraicspaceswithanactionofGm}. For a positive integer $n$, the natural map $X^{\lambda,+}\to X^{\lambda^n,+}$ is an isomorphism (this can be seen as a consequence of \cite[Thm. 5.1.4]{Bu-HalpernLeistner-IbanezNunez-Kinjo-IntrinsicDTtheoryI}), and thus $X^{\lambda,+}$ is defined also for rational cocharacters $\mathbb{X}_*(T)_\bQ$. 
\end{para}

\begin{para}[Stacks of rational filtrations and graded points]
    For given $\lambda\in \mathbb{X}_*(T)_\bQ$, the group $P_\lambda$ acts on $X^{\lambda,+}$ and $L_\lambda$ acts on $X^{\lambda,0}$. Let $C\subset \mathbb{X}_*(T)_\bQ$ be a complete set of representatives for $\mathbb{X}_*(T)_\bQ/W$.  We define the stack of \emph{rational filtrations} of $X/G$ to be the disjoint union
    \[\Filt_\Q(X/G)=\bigsqcup_{\lambda\in C} X^{\lambda,+}/P_\lambda.\]
    Similarly, the stack of \emph{rational graded points} of $X/G$ is
    \[\Grad_\Q(X/G)=\bigsqcup_{\lambda\in C} X^{\lambda,0}/L_\lambda.\]
\end{para}

\begin{remark}[Intrinsic definition]
     The stacks of rational filtrations and graded points are intrinsic to $X/G$ and do not depend on the presentation as quotient stack nor on the choice of complete set of representatives $C$. Beyond quotient stacks, they can be defined as certain colimit of mapping stacks from $\A^1_k/\G_{m,k}$ and $\mathrm{B}\G_{m,k}$. See \cite[Theorem 1.8.4]{_HalpernLeistner_Onthestructureofinstabilityinmodulitheory} and \cite[Section 2.2]{ibanez-thesis}.
\end{remark}
   
\begin{para}[Evaluation map and set of filtrations]
    The forgetful map $X^{\lambda,+}\to X$ defined by precomposition along $\{1\}\to \A^1_k$ is equivariant with respect to the homomorphism $P_\lambda\to G$. Thus, there is a natural forgetful map
    \[\ev\colon \Filt_\Q(X/G)\to X/G.\]
    
    The set of filtrations $\qfilt(X/G,x)$ is identified with the set of $k$-points of the fibre of $\ev$ at $x\colon \Spec k\to X/G$, and thus it is independent of the presentation of $X/G$ as a quotient stack. 
\end{para}

\begin{para}[Associated graded map]
    The map $\gr\colon X^{\lambda,+}\to X^{\lambda,0}$ defined by precomposition along $\{0\}\to \A^1_k$ is equivariant with respect to the homomorphism $P_\lambda\to L_\lambda$ and thus it defines an \emph{associated graded} morphism
    \[\gr\colon \Filt_\Q(X/G)\to \Grad_\Q(X/G).\]
\end{para}

\begin{para}[The component lattice]
    We denote $\abs{\CL_\Q(X/G)}$ the set of connected components of $\Grad_\Q(X/G)$. This set can be endowed with additional combinatorial structure giving rise to what is called the \emph{component lattice} of $X/G$. Component lattices of general stacks are studied by Bu, Halpern-Leistner, Kinjo and the second named author \cite{Bu-HalpernLeistner-IbanezNunez-Kinjo-IntrinsicDTtheoryI} with a view towards applications in enumerative geometry. We will not need the full structure of the component lattice in this note.
\end{para}

\begin{para}[Notations]
    Explicitly, an element $\alpha \in \abs{\CL_\Q(X/G)}$ is a pair $\alpha=(\lambda_\alpha, X_\alpha)$ with $\lambda_\alpha\in C$ and $X_\alpha$ a closed and open $L_{\lambda_\alpha}$-equivariant subscheme of $X^{\lambda_\alpha,0}$ such that $X_\alpha/L_{\lambda_\alpha}$ is connected (that is, the data of a connected component of the stack $X^{\lambda_\alpha,0}/L_{\lambda_\alpha}$). We denote $L_\alpha\coloneqq L_{\lambda_\alpha}$ and $P_\alpha\coloneqq P_{\lambda_\alpha}$.

    The associated graded morphism $\gr\colon \Filt_\Q(X/G)\to \Grad_\Q(X/G)$ induces a bijection on connected components.
    Therefore, for every $\alpha\in \abs{\CL_\Q(X/G)}$ the preimage of $X_\alpha/L_\alpha$ along $\gr$ is a connected component of $\Filt_\Q(X/G)$, of the form $X^+_\alpha/P_\alpha$, where $X^+_\alpha$ is a closed and open $P_\alpha$-equivariant subscheme of $X^{\lambda_\alpha,+}$.
\end{para}

\begin{para}[Hilbert-Mumford weight and the component lattice]
    Let $\alpha\in \abs{\CL_\Q(X/G)}$ and take a point $x\in X_\alpha(k)$. The number $\mathrm{m}(x,\lambda_\alpha)$ does not depend on the choice of point $x\in X_\alpha(k)$, so we may define $\mathrm{m}(\alpha)\coloneqq \mathrm{m}(x,\lambda_\alpha)$. We will also use the notation $\norm{\alpha}\coloneqq \norm{\lambda_\alpha}$.
\end{para}

\subsection{Kirwan stratification}

\begin{para}[Stability function]
    We define the \emph{stability function} $M\colon X(k)\to \R_{\geq 0}\cup \{\infty\}$ by
    \[M(x)=\sup \left\{ \norm{\lambda} \left\vert \ \lambda\in\qfilt(X/G,x),\ \mathrm{m}(x,\lambda)=\norm{\lambda}^2\right.\right\},\]
    for $x\in X(k)$.
    In principle $M(x)$ could be $\infty$, but we will see that this is never the case in our framework.

    Note that a point $x\in X(k)$ is semistable if and only if $M(x)=0$.
\end{para}

\begin{theorem}[Kirwan stratification] \label{thm: kirwan stratification} The following statements hold:
    \begin{enumerate}[topsep=\parasep/2]
    \item\label{item-definition-Xalpha+ss} For every $\alpha \in \abs{\CL_\Q(X/G)}$, there is a unique $P_\alpha$-equivariant open subscheme $X_\alpha^{+\ss}$ of $X^+_\alpha$ such that, for every $x\in X^+_\alpha(k)$, the point $x$ is in $X_\alpha^{+\ss}$ if and only if $M(\mathrm{ev}(x))=\norm{\alpha}$ and $\mathrm{m}(\alpha)=\norm{\alpha}^2$.
    \item The maps $\mathrm{ev}\colon X_\alpha^{+\ss}/P_\alpha=\left(G\times^{P_\alpha}X_\alpha^{+\ss}\right)/G\to X/G$ are pairwise disjoint locally closed immersions and jointly surjective. 
    \item For every $c\in \R$, the set 
\[\bigsqcup_{\substack{\alpha\in \abs{\CL_\Q(X/G)} \\ \norm{\alpha}\geq c}}\abs{G\times^{P_\alpha} X_\alpha^{+\ss}}\] 
    is closed in $\abs{X}$.
    \item There is a unique $L_\alpha$-equivariant open subscheme $X_\alpha^\ss$ of $X_\alpha$ such that $\gr_\alpha^{-1}(X_\alpha^\ss)=X_\alpha^{+\ss}$.
    \item The quotient stack $X_{\alpha}^{\ss}/L_\alpha$ admits a good moduli space $X_{\alpha}^{\ss} \git L_\alpha$ which is projective over the scheme $S^{\alpha,0}\git L_\alpha$. In particular, $X_{\alpha}^{\ss} \git L_\alpha$ is projective over $S\git G$.
    \end{enumerate}
\end{theorem}
\begin{proof}
    In the language of \cite{_HalpernLeistner_Onthestructureofinstabilityinmodulitheory} and \cite{ibanez-thesis}, the rational line bundle $\cL$ defines a \emph{linear form $\ell$ on graded points} of $X/G$, 
    and the $W$-invariant rational inner product on $\mathbb{X}_*(T)_\bQ$ defines a \emph{norm $q$ on graded points} of $X/G$. 
    
    The first four statements in the theorem mean precisely that $\ell$ and $q$ define a $\Theta$-stratification of $X/G$, which is a particular case of \cite[Thm. 2.6.4]{ibanez-thesis}.

    The fifth statement is true if $\lambda_\alpha=0$ by \cref{para-good-moduli-ss-locus}. For general $\alpha$ it reduces to \cref{para-good-moduli-ss-locus} by \cref{thm: description of centres} below. 
    The last sentence in the fifth statement follows from the fact that the morphism $S^{\alpha,0}\git L_\alpha \to S \git G$ is finite \cite[Thm. in page 234]{luna_orbite}.
\end{proof}

\begin{para}[Harder--Narasimhan filtration]
\label{para-HN filtration}
Every point $x\in X(k)$ belongs to a unique stratum, and thus there is a unique $\alpha\in |\CL_\bQ(X/G)|$, a unique point $y\in (X_\alpha^{+\ss}/P_\alpha)(k)$ and an isomorphism $y\simeq x$ in the groupoid $(X/G)(k)$. The data defines an element $\mathrm{HN}(x)\in \qfilt(X/G,x)$ called the \emph{Harder--Narasimhan filtration} of $x$. It is the unique maximizer of the norm function $\norm{\lambda}$ on $\qfilt(X/G,x)$ subject to the condition that $\mathrm{m}(x,\lambda)=\norm{\lambda}^2$. 
\end{para}

\begin{lemma}
For every $\lambda\in \mathbb{X}_*(T)_\bQ$ there is a unique rational cocharacter $\lambda^\vee\in \mathbb{X}^*(L_\lambda)_\bQ$ such that, for every $\eta\in \mathbb{X}_*(T)_\bQ$, we have
\[\langle \eta, \lambda^\vee \rangle=(\eta,\lambda).\]
\end{lemma}
Here $\langle -,- \rangle$ denotes the duality paring $\mathbb{X}_*(T)_\bQ\times \mathbb{X}^*(T)_\bQ\to \Q$.
\begin{proof}
Since the inner product $(-,-)$ is Weyl-invariant, the map $(-,\lambda)\colon \mathbb{X}_*(T)_\bQ\to \Q$ is invariant under the action of the Weyl group $W(T,L_\lambda)$ of $T$ inside $L_\lambda$. Therefore, by \Cref{theorem: rational characters and Weyl group}, there is a unique rational cocharacter $\lambda^\vee$ such that $\lambda^\vee=(-,\lambda)$.
\end{proof} 

\begin{para}[Shifted linearization]
    Let $\alpha\in \abs{\CL_\Q(X/G)}$ and denote $\alpha^\vee=\lambda_\alpha^\vee\in \mathbb{X}^*(L_\alpha)_\bQ$. The rational character $\alpha^\vee$ of $L_\alpha$ defines a rational line bundle $\cO_{\mathrm{B}L_\alpha}(\alpha^\vee)$ on $\mathrm{B}L_\alpha$. Let $\cO_{X_\alpha/L_\alpha}(\alpha^\vee)$ be the pullback of $\cO_{\mathrm{B}L_\alpha}(\alpha^\vee)$ along $X_\alpha/L_\alpha\to \mathrm{B}L_\alpha$. The \emph{shifted linearization} on $X_\alpha/L_\alpha$ is the rational line bundle $\cL(\alpha^\vee)=\cL\vert_{X_\alpha/L_\alpha}\otimes \cO_{X_\alpha/L_\alpha}(\alpha^\vee)$ on $X_\alpha/L_\alpha$.
\end{para}

We emphasize that the existence of the shifted linearization, and thus the following theorem, crucially depend on \Cref{theorem: main intro thm}.

\begin{theorem}[Description of the centres as semistable loci] \label{thm: description of centres}
The semistable locus of $X_\alpha$ with respect to the shifted linearization $\cL(\alpha^\vee)$ is precisely $X_\alpha^\ss$.
\end{theorem}
\begin{proof}



The Recognition Theorem \cite[Thm. 5.4.4]{_HalpernLeistner_Onthestructureofinstabilityinmodulitheory} states that a point $x\in X_\alpha(k)$ is in $X_\alpha^\ss$ if and only if $\lambda_\alpha$ maximizes the function $u(\lambda)= \frac{\mathrm{m}(x,\lambda)}{\norm{\lambda}}$ on $0\neq \lambda\in \qfilt(X_\alpha/L_\alpha,x)\setminus \{0\}$, or, equivalently, if $\lambda_\alpha$ maximizes $\norm{\lambda}$ for $\lambda\in\qfilt(X_\alpha/L_\alpha,x)$ subject to the equality $\mathrm{m}(x,\lambda)=\norm{\lambda}^2$. Indeed, we can scale $\lambda$ to make the equality hold, and in that case $u(\lambda)=\norm{\lambda}$. Yet another equivalent condition is that $\lambda_\alpha$ maximizes the function $f(\lambda)=\mathrm{m}(x,\lambda)-\frac{1}{2}\norm{\lambda}^2$ on $\qfilt(X_\alpha/L_\alpha,x)$. Indeed, evaluating the derivative of $h(t)=f(t\lambda)$ at $t=1$ we see that any $\lambda$ maximizing $f$ must satisfy $\mathrm{m}(x,\lambda)=\norm{\lambda}^2$, while for the elements $\lambda$ satisfying this equality we have $f(\lambda)=\dfrac{1}{2}\norm{\lambda}^2$.

 Now, since $\lambda_\alpha$ is central in $L_\alpha$ and fixes $x$, for any other $\lambda\in\qfilt(X_\alpha/L_\alpha,x)$ and any $t\in \bQ_{\geq 0}$, we have $\lambda_\alpha+t\lambda\in \qfilt(X_\alpha/L_\alpha,x)$. Thus $x$ is in $X_\alpha^\ss$ if and only if for every $\lambda\in\qfilt(X_\alpha/L_\alpha,x)$ the function $v(t)=f(\lambda_\alpha+t\lambda)$ (defined for $t\geq 0$) has a maximum at $t=0$. Since the function $v(t)=\mathrm{m}(x,\lambda_\alpha)+t\mathrm{m}(x,\lambda)-\frac{1}{2}\norm{\lambda_\alpha}^2-\frac{1}{2}t^2\norm{\lambda}^2-(\lambda,\lambda_\alpha)t$ is concave and differentiable, this happens precisely if $v'(0)\leq 0$. But $v'(0)= m(x,\lambda)-\langle\lambda,\lambda_\alpha\rangle$, which is the Hilbert-Mumford weight for the shifted linearization $\cL(\alpha^\vee)$. Thus the condition is equivalent to semistability of $x$ with respect to $\cL(\alpha^\vee)$.
\end{proof}

The argument above goes through for general stacks, giving the \emph{Linear Recognition Theorem} in \cite{HalpernLeistner-IbanezNunez-variation-instability}.

\begin{remark}
    In the case when $S = \Spec(k)$, the statement of \Cref{thm: description of centres} is claimed in \cite[Remark 2.21]{_Kirwan_Cohomologyofquotientsinsymplecticandalgebraicgeometry}. In that argument, it is assumed that there is a rational character of a Levi subgroup $L_{\lambda}$ that is dual to a given cocharacter $\lambda$, but no proof is provided. There seems to be no complete and detailed argument of this description of the stratification in the literature. Our arguments provide a proof of this.
\end{remark}

\section{Semistability for principal bundles} \label{section: semistability of G bundles}
\addtocounter{subsection}{1}
In this section, we use the results developed in Sections \cref{section: weyl groups and cocharacters} and \cref{section: kirwan stratification} to study the notion of semistability of $G$-bundles for a disconnected reductive group $G$.

\begin{para}[Setup for the section]
We keep the assumption that the ground field $k$ is algebraically closed of characteristic $0$. Fix a reductive group $G$ and a smooth connected projective curve $C$ over $k$. We fix once and for all a maximal torus $T \subset G$, and denote by $W_G$ the Weyl group of $G$. We denote by $\Bun_G=\mathrm{Map}(C,\mathrm{B}G)$ the stack of $G$-bundles on the curve $C$, which is a smooth algebraic stack over $k$ \cite[Proposition 1]{heinloth-uniformization}.
\end{para}

\subsection{Rational degree}

\begin{para}[The stack $\Bun_{\pi_0(G)}$]
    Let $G^{\circ} \subset G$ denote the neutral component of $G$, and let $\pi_0(G) = G/G^{\circ}$ be the group of connected components of $G$. Standard deformation theory as in \cite[Proposition 1]{heinloth-uniformization} shows that the stack $\Bun_{\pi_0(G)}$ is an \'etale Deligne-Mumford stack, and its topological space $|\Bun_{\pi_0(G)}|$ is a discrete union of points corresponding to isomorphism classes of $\pi_0(G)$-bundles on $C$. Thus $\Bun_{\pi_0(G)}$ is a disjoint union of stacks of the form $\mathrm{B}H$ with $H$ a finite discrete group.

    When $k=\bC$, we can write $\Bun_{\pi_0(G)}=\Hom(\pi_1(C),\pi_0(G))/\pi_0(G)$ explicitly as a quotient stack, where $\Hom(\pi_1(C),\pi_0(G))$ is regarded as a disjoint union of copies of $\Spec \bC$ and $\pi_0(G)$ acts by conjugation.
\end{para}

\begin{para}[Components of $\Bun_{\pi_0(G)}$]
    \label{para: components of pi0(G)} 
    Note that each $\pi_0(G)$-bundle $E$ on $C$ can be written as the extension of structure group $E'(\pi_0(G))$ of a connected $F$-bundle $E'$ over $C$, where $F \subset \pi_0(G)$ is a subgroup of $F$ that is uniquely determined up to conjugation (one can take $F$ the subgroup that preserves a chosen connected component of $E$). This gives a map $\vert \Bun_{\pi_0(G)} \vert\to \mathrm{Conj}(\pi_0(G))$, where $\mathrm{Conj}(\pi_0(G))$ is the set of conjugacy classes of subgroups of $\pi_0(G)$. When working over $\bC$, if a $\pi_0(G)$-bundle $E$ is given by a homomorphism $\phi\colon \pi_1(C) \to \pi_0(G)$, then the corresponding subgroup $F\subset \pi_0(G)$ is the image of $\phi$, which depends only on $E$ up to conjugation. 

    We also have a finite \'etale and surjective morphism
    \[ \bigsqcup_{[F] \in \mathrm{Conj}(\pi_0(G))} \Bun_{F}^{\circ} \to \Bun_{\pi_0(G)},\]
    where $\Bun_{F}^{\circ} \subset \Bun_F$ is the open and closed substack parametrizing $F$-bundles over $C$ that are connected.
\end{para}

\begin{para}[Decomposition of $\Bun_G$ by conjugacy classes of subgroups in {$\pi_0(G)$}]
    The quotient morphism $G \twoheadrightarrow \pi_0(G)$ induces a morphism of stacks $\Bun_G \to \Bun_{\pi_0(G)}$ given by taking a $G$-bundle $E$ to the associated $\pi_0(G)$-bundle $E(\pi_0(G))$. Given a subgroup $F \subset \pi_0(G)$, we denote by $\Bun_G^{F}$ the fiber product $\Bun_G \times_{\Bun_{\pi_0(G)}} \Bun_F^{\circ}$, where the morphism $\Bun_F^{\circ} \to \Bun_{\pi_0(G)}$ is as in \cref{para: components of pi0(G)}. By construction, the induced morphism $\Bun_G^F \to \Bun_G$ is finite and \'etale. If we denote by $G^F \subset G$ the preimage of $F$ under the quotient morphism $G \to \pi_0(G)$, then $\Bun_G^F$ is an open and closed substack of the stack $\Bun_{G^F}$, namely the open substack parametrizing $G^F$-bundles such that associated $F$-bundle on $C$ is connected.
\end{para}

\begin{para}[Set of rational degrees]
    Let $S_G$ be the set of pairs $(F,d)$, where $F \subset \pi_0(G)$ is a subgroup and $d \in \mathbb{X}^*(G^F)_\bQ^{\vee}$. The group $G(k)$ acts on $S_G$ by conjugation. We define the \emph{set $\pi_1(G)_\bQ$ of rational degrees} for $G$ to be the set quotient
        \[\pi_1(G)_\bQ\coloneqq S_G/G(k).\]
    Note that $G^{\circ}(k)$ acts trivially, so we get an action of $\pi_0(G)$ on $S_G$ and $\pi_1(G)_\bQ= S_G/\pi_0(G)$. For $(F,d)\in S$, we denote $[F,d]$ the corresponding element in $\pi_1(G)_\bQ$.
\end{para}

\begin{para}[Rational degree of a $G$-bundle]
    \label{para: rational degree of G-bundle}
    Given an algebraically closed field extension $K \supset k$ and a $K$-point $E \in \Bun_G(K)$ corresponding to a $G$-bundle on $C_K$, we define its \emph{rational degree} $[F_E, d_E] \in \pi_1(G)_{\mathbb{Q}}$ as follows. First, we take a subgroup $F\subset \pi_0(G)$ and a lift $\widetilde{E} \in \Bun_G^F(K)$ of $E$, thought of as a $G^F$-bundle. We define the linear functional $d$ on $\mathbb{X}^*(G^F)_\bQ$ which sends a rational character $\chi \in \mathbb{X}^*(G^F)_\bQ$ to the degree of the associated rational line bundle $\widetilde{E}(\chi)$ on $C_K$. The pair $(F,d)\in S_G$ is well-defined up to the action of $\pi_0(G)$, and we let $[F_E,d_E]=[F,d]\in \pi_1(G)_{\mathbb{Q}}$. 
    
    The rational degree $[F_E,d_E]$ is preserved under base change to any algebraically closed field extension of $K$. Also, since the image of each $\Bun_G^F \to \Bun_G$ is open and closed, and since the degree of rational line bundles is locally constant in families, the assignment
    \[ \vert\Bun_G\vert \to \pi_1(G)_{\mathbb{Q}}, \; \; E \mapsto [F_E, d_E]\]
    is a locally constant function on the underlying topological space $|\Bun_G|$ of $\Bun_G$. For any given $[F,d] \in \pi_1(G)_{\mathbb{Q}}$, we denote by $\Bun_G^{[F,d]} \subset \Bun_G$ the open and closed substack of $\Bun_G$ parametrizing $G$-bundles of rational degree $[F,d]$.
\end{para}

\begin{remark}
    While we could have defined the rational degree of $E$ simply as an element of $\mathbb{X}^*(G)_\bQ^\vee$, without taking $F$ into account, this alternative notion would not contain enough information for our purposes. 
    \cref{example: subgroups of pi0 are necessary} illustrates this defect.
\end{remark}

\begin{example}[Connected reductive groups] \label{example: rational degrees for connected groups}
    If $G$ is a connected reductive group, then $\pi_1(G)_{\mathbb{Q}} \cong \mathbb{X}^*(G)_\bQ^{\vee} = \mathbb{X}_*(Z(G))_\bQ= \mathbb{X}_*(T)_\bQ^{W_G}$ is a vector space. Given a $G$-bundle $E \in \Bun_G(K)$, the corresponding degree $d_E \in \pi_1(G)_{\mathbb{Q}}$ is the linear functional that sends a rational character $\chi \in \mathbb{X}^*(G)_\bQ$ to the degree $\deg(E(\chi))$.
\end{example}

\begin{para}[Rational degree under change of group] 
    \label{para: rational degree change group}
    Suppose that we are given another smooth linearly reductive affine group $H$ and a homomorphism $f: G \to H$. We get an induced map $S_f\colon S_G\to S_H$ as follows. For $(F,d) \in S_G$, let $F' = f(F) \subset \pi_0(H)$ be the image of $F$ under the map  $\pi_0(G) \to \pi_0(H)$ on connected components. The induced homomorphism $f_F\colon G^F\to H^{F'}$ gives a map $\mathbb{X}^*(G^F)_\bQ^{\vee} \to \mathbb{X}^*(H^{F'})_\bQ^{\vee}$ sending $d$ to an element $d' \in \mathbb{X}^*(H^{F'})_\bQ^{\vee}$. We let $S_f(F,d) = (F',d') \in S_H$. Since $S_f$ is equivariant with respect to the homomorphism $G(k)\to H(k)$, it descends to a map $\pi_1(f): \pi_1(G)_{\mathbb{Q}} \to \pi_1(H)_{\mathbb{Q}}$.

    There is an induced morphism of stacks $f_*: \Bun_G \to \Bun_H$ which sends a $G$-bundle $E$ to its associated $H$-bundle $E(H)$. Given a degree $[F,d] \in \pi_1(G)_{\mathbb{Q}}$, the image of the open and closed substack $\Bun_G^{[F,d]} \subset \Bun_G$ under $f_*: \Bun_G \to \Bun_H$ lies inside $\Bun_H^{\pi_1(f)([F,d])}$.
\end{para}

\begin{definition}[Adapted degree]
    \label{definition: adapted degree}
    Let $f\colon G\to H$ be a homomorphism as in \cref{para: rational degree change group} and let $[F,d]\in \pi_1(G)_\bQ$ be a degree. We denote $f_F\colon G^F\to H^{f(F)}$ the induced map and we identify $\mathbb{X}^*(G^F)_\bQ^{\vee}\simeq \mathbb{X}_*(Z(G^F))_\bQ$ by \Cref{theorem: rational characters and Weyl group}.  Let $d'$ denote the image of $d$ under the induced map $\mathbb{X}_*(Z(G^F))_\bQ\to \mathbb{X}_*(H^{f(F)})_\bQ$ on cocharacters. We say that the degree $[F,d]$ is \emph{$f$-adapted} if $d'$ lies in the subspace $\mathbb{X}_*(Z(H^{f(F)}))_\bQ\subset \mathbb{X}_*(H^{f(F)})_\bQ$ of central cocharacters. This condition does not depend on the choice of representative $(F,d)\in S_G$ of $[F,d]$. 
\end{definition}

\begin{para}
    \label{remark: enough to consider centre of connected component}
    Note that, since the induced map $\pi_0(G^F)\to \pi_0(H^{f(F)})$ on components is surjective, the condition is equivalent to $d'$ lying in the bigger subspace $\mathbb{X}_*(Z(H^{\circ}))_\bQ\subset \mathbb{X}_*(H^{f(F)})_\bQ$ of cocharacters that are central in $H^\circ$.
\end{para}

\begin{para}
    \label{remark: assume pi_0(f) surjective}
    Given $(F,d)\in S_G$, the degree $[F,d]\in \pi_1(G)_\bQ$ is $f$-adapted if and only if the lift $[F,d]\in \pi_1(G^F)_\bQ$ is $f_F$-adapted. 
\end{para}

\begin{para}[In terms of Weyl groups]
    \label{para: in terms of Weyl groups}
    Given a homomorphism $f\colon G\to H$ as in \cref{para: rational degree change group}, we can express the condition of $f$-adaptedness in terms of Weyl groups. We fix a maximal torus $T' \subset H$ containing $f(T)$, thus getting a homomorphism of vector spaces $\mathbb{X}_*(T)_\bQ \to \mathbb{X}_*(T')_\bQ$. Let $[F,d]\in\pi_1(G)_\bQ$ be a degree. By \cref{remark: assume pi_0(f) surjective}, after replacing $f$ by $f_F$ we may assume that $F=\pi_0(G)$ and that $\pi_0(G)\to \pi_0(H)$ is surjective.

    By \Cref{theorem: cocharacters and Weyl group,theorem: rational characters and Weyl group}, we can identify $\mathbb{X}^*(G)_\bQ^{\vee} =\mathbb{X}_*(Z(G))_\bQ= \mathbb{X}_*(T)_\bQ^{W_{G}}$ as a subspace of $\mathbb{X}_*(T)_\bQ$ and $\mathbb{X}_*(Z(H))_\bQ=\mathbb{X}_*(T')_\bQ^{W_H}$ as a subspace of $\mathbb{X}_*(T')_\bQ$. Thus the degree $[F,d]$ is $f$-adapted precisely if the image $d'$ of $d$ under $\mathbb{X}_*(T)_\bQ\to \mathbb{X}_*(T')_\bQ$ lies in the subspace $\mathbb{X}_*(T')_\bQ^{W_H}\subset \mathbb{X}_*(T')_\bQ$. By \cref{remark: enough to consider centre of connected component}, this condition is equivalent to $d'$ lying in the subspace  $\mathbb{X}_*(T')_\bQ^{W_{H}^{\circ}} \subset \mathbb{X}_*(T')_\bQ$ of invariants under the Weyl group $W_{H}^{\circ} \subset W_H$ of the neutral component $H^{\circ} \subset H$.
\end{para}

\begin{para}
    \label{para: decomposition of map on rational degrees}
    Using the notation of the previous paragraph, we have a decomposition of $\mathbb{X}^*(G)_\bQ^{\vee}\to \mathbb{X}^*(H)_\bQ^\vee$ as
    \[ \mathbb{X}^*(G)_\bQ^{\vee}= \mathbb{X}_*(T)_\bQ^{W_G} \hookrightarrow \mathbb{X}_*(T)_\bQ \to \mathbb{X}_*(T')_\bQ \twoheadrightarrow \mathbb{X}_*(T')_\bQ^{W_{H}} = \mathbb{X}^*(H)_\bQ^{\vee},\]
    where $\mathbb{X}_*(T')_\bQ \twoheadrightarrow \mathbb{X}_*(T')_\bQ^{W_{H}}$ is the averaging operator. Thus $d$ being $f$-adapted means that there is no need to project $d'$ using $\mathbb{X}_*(T')_\bQ \twoheadrightarrow \mathbb{X}_*(T')_\bQ^{W_{H}}=\mathbb{X}^*(H)_\bQ^\vee$ to obtain the corresponding degree for $H$.
\end{para}

\begin{example}[Direct sum homomorphism] \label{example: rational degrees for direct sum homomorphism}
    Fix two integers $n,m>0$. Set $G = \GL_n \times \GL_m$ and $H = \GL_{n+m}$. Define $f: G \to H$ to be the direct sum homomorphism, which sends a pair of matrices $(A,B)$ to the block diagonal matrix 
    \[\begin{bmatrix} A & 0 \\ 0 & B\end{bmatrix}.\] 
    We have a canonical identification $\pi(G)_{\mathbb{Q}} = \mathbb{X}^*(G)_\bQ^{\vee} = \mathbb{Q}^{2}$ via the determinant characters, and similarly $\pi_1(H)_{\mathbb{Q}}= \mathbb{Q}$. 
        
    A $G$-bundle $E$ on $C$ corresponds to a pair of vector bundles $E= (E_1,E_2)$, where $E_1$ has rank $n$ and $E_2$ has rank $m$. The degree $d_E$ is given by the pair $(\deg(E_1), \deg(E_2)) \in \mathbb{Q}^{\oplus 2} = \pi_1(G)_{\mathbb{Q}}$ of degrees of the vector bundles. Similarly, $H$-bundles $E'$ correspond to rank $(n+m)$ vector bundles, and the corresponding degree $d_{E'} \in \mathbb{Q}$ is the degree of the vector bundle. The induced morphism of stacks $f_*: \Bun_G \to \Bun_H$ sends a pair $(E_1,E_2)$ of vector bundles to their direct sum $E_1 \oplus E_2$. The homomorphism $\pi_1(G)_{\mathbb{Q}} = \mathbb{Q}^{\oplus 2} \to \mathbb{Q} = \pi_1(H)_{\mathbb{Q}}$ sends a pair $(a,b)$ of rational numbers to their sum $a+b$.

    Fix the diagonal maximal tori $T \subset G$ and $T'\simeq T \subset H$. We use the standard basis to write $\mathbb{X}_*(T')_\bQ = \mathbb{Q}^{n+m} = \mathbb{Q}^{n} \oplus \mathbb{Q}^{m}$. The corresponding composition
    \[ \mathbb{Q}^{2} = \pi_1(G)_{\mathbb{Q}} \hookrightarrow \mathbb{X}_*(T)_\bQ = \mathbb{X}_*(T')_\bQ = \mathbb{Q}^{n} \oplus \mathbb{Q}^{m}\]
    sends a pair $(a,b)$ to the pair of vectors $(v_1,v_2) \in \mathbb{Q}^{n} \oplus \mathbb{Q}^{m}$ given by
    \[ v_1= \left( \frac{a}{n}, \frac{a}{n}, \ldots , \frac{a}{n}\right) \; \; \text{ and } \; \; \; v_2= \left( \frac{b}{m}, \frac{b}{m}, \ldots, \frac{b}{m}\right).\]
    Hence, a degree $(a,b)$ is $f$-adapted if and only if $\frac{a}{n} = \frac{b}{m}$. In terms of bundles, the degree of a pair of vector bundles $(E_1, E_2)$ is $f$-adapted if and only if the slopes of the vector bundles $E_1$ and $E_2$ are equal.
    \end{example}

\begin{para}[Adapted degree in the connected case]
When $G$ is connected, \textcite{Schieder_2014} refers to the map $s\colon\pi_1(G)_\bQ \cong \bX_*(T)_\bQ^{W_G}\hookrightarrow \bX_*(T)$ as the \emph{slope} map. If $f\colon G\to H$ is a homomorphism as above, $H$ is also connected, and $T'\subset H$ is a maximal torus containing $f(T)$, then $f$ induces natural maps $f\colon \pi_1(G)_\bQ\to \pi_1(H)_\bQ$ and $f\colon \bX_*(T)\to \bX_*(T')$. Using this language, a degree $d\in \pi_1(G)_\bQ$ is $f$-adapted precisely when $s(f(d))=f(s(d))$.
We thank Dragoș Frățilă for pointing this out.

\end{para}

\subsection{Hilbert-Mumford semistability}

\begin{para}[Determinant of cohomology line bundle]
    Let $\Ad: G \to \GL(\mathfrak{g})$ denote the adjoint representation of $G$ on its Lie algebra $\mathfrak{g}$. Let $E_{\mathrm{univ}}$ denote the universal $G$-bundle on $\Bun_G \times C$, and let $E_{\mathrm{univ}}(\mathfrak{g})$ denote the associated adjoint vector bundle. Following \textcite{_Heinloth_HilbertMumfordstabilityonalgebraicstacksandapplicationstoGbundlesoncurves}, we define the line bundle $\mathcal L_{\det}$ on the stack $\Bun_G$ by
    \[\mathcal L_{\det} := \det(R\pi_*E_{\mathrm{univ}}(\mathfrak{g})),\]
     where $\pi: \Bun_G \times C \to \Bun_G$ denotes the first projection, and we take the determinant of the perfect complex $R\pi_*E_{\mathrm{univ}}(\mathfrak{g})$ on $\Bun_G$.
\end{para}

Now we may follow \cite[\S1.F]{_Heinloth_HilbertMumfordstabilityonalgebraicstacksandapplicationstoGbundlesoncurves} to define a notion of semistability for $G$-bundles using the Hilbert-Mumford weight for the line bundle $\mathcal L_{\det}$.

\begin{definition}[Semistability of $G$-bundles] \label{defn: semistability bundles}
    Let $K \supset k$ be an algebraically closed field extension, and let $E \in \Bun_G(K)$ be a $G$-bundle on $C_K$. We say that $E$ is (Hilbert-Mumford) \emph{semistable} if for all morphisms $\varphi\colon \mathbb{A}^1_K/\mathbb{G}_{m,K} \to \Bun_G$ such that $\varphi(1) \cong E$, we have $\wt(\varphi^*(\mathcal L_{\det})|_0) \leq 0$.
\end{definition}

\begin{para}[Filtrations and graded points]
    \label{para: filtrations and graded points}
    In the framework of \cref{defn: semistability bundles}, a morphism $\varphi\colon \mathbb{A}^1_K/\mathbb{G}_{m,K} \to \Bun_G$ with $\varphi(1) \cong E$ is called a \emph{filtration} of $E$, following the language of \textcite{_HalpernLeistner_Onthestructureofinstabilityinmodulitheory}. Specifying a filtration of $E$ is equivalent to specifying a one-parameter subgroup $\lambda\in \mathbb{X}_*(T)$ together with a reduction of structure group $E_P$ of $E$ to $P=P_G(\lambda)$, since
    \[\mathrm{Map}(\bA^1_K/\mathbb{G}_{m,K},\mathrm{Map}(C,\mathrm{B}G))=\mathrm{Map}(C,\mathrm{Map}(\bA^1_K/\mathbb{G}_{m,K},\mathrm{B}G))\]
    and 
    \[\mathrm{Map}(\bA^1_K/\mathbb{G}_{m,K},\mathrm{B}G)=\mathrm{Filt}(\mathrm{B}G)=\bigsqcup_{\lambda\in \mathbb{X}_*(T)/W}\mathrm{B}P_G(\lambda),\]
    by \cite[Theorem~1.4.8]{_HalpernLeistner_Onthestructureofinstabilityinmodulitheory}.

    A \emph{graded point} of $\Bun_G$ is a morphism $\mathrm{B}\mathbb{G}_{m,K}\to \Bun_G$. Similarly to the case of filtrations, a graded point $\varphi: \mathrm{B}\mathbb{G}_{m,K} \to \Bun_G$ is determined by the data of a cocharacter $\lambda \in \mathbb{X}_*(T)$ and an $L_G(\lambda)$-bundle.

    If $\varphi\colon \mathbb{A}^1_K/\mathbb{G}_{m,K} \to \Bun_G$ is a filtration of $E$, then the restriction $\varphi\vert_{\mathrm{B}\mathbb{G}_{m,K}}\colon \mathrm{B}\mathbb{G}_{m,K}\to \Bun_G$ is a graded point, called the \emph{associated graded point} of the filtration $\phi$. It corresponds to the $L_G(\lambda)$-bundle $E_P(L_G(\lambda))$ associated to $E_P$ under the quotient map $P_G(\lambda) \to L_G(\lambda)$.
\end{para}
\begin{proposition}[Semistability under change of group of components] \label{prop: behavior stability under change of group of components}
    Let $K \supset k$ be an algebraically closed field extension and let $E \in \Bun_G(K)$ be a $G$-bundle on $C_K$. Let $F \subset \pi_0(G)$ be a subgroup, and consider the preimage $G^F \subset G$. Suppose that there is a reduction of structure group of $E$ to a $G^F$-bundle $\widetilde{E}$. Then, $E$ is semistable if and only if $\widetilde{E}$ is semistable.
\end{proposition}
\begin{proof}
    Let $i: G^F \hookrightarrow G$ denote the inclusion. Since the morphism $i_*: \Bun_{G^F} \to \Bun_G$ is finite, the valuative criterion for properness implies that postcomposing with $i_*$ induces a bijection between the set of filtrations of $\widetilde{E}$ and the set of filtrations of $E$. Since the line bundle $\cL_{\det}$ on $\Bun_{G^F}$ is the pullback $i^*(\cL_{\det})$ of the line bundle $\cL_{\det}$ on $\Bun_G$, it follows that the weight requirements in \Cref{defn: semistability bundles} for the semistability of $\widetilde{E}$ are identical to the weight requirements for the semistability of $E$.
\end{proof}

The rest of this subsection is dedicated to explaining a reformulation of \Cref{defn: semistability bundles} in terms of parabolic reductions, following \textcite{_Heinloth_HilbertMumfordstabilityonalgebraicstacksandapplicationstoGbundlesoncurves} and F. H. and Zhang \cite{herrero2023meromorphic}. We first need some setup. 

\begin{para}[Trace pairing]
    The adjoint representation $\mathfrak{g}$ induces a bilinear trace form on the group $\mathbb{X}_*(T)_\bQ$ given as follows. Let $\mathfrak{g} = \bigoplus_{ \alpha \in \mathbb{X}^*(T)} \mathfrak{g}_{\alpha}= \bigoplus_{ \alpha \in \Phi} \mathfrak{g}_{\alpha}$ be the decomposition of $\mathfrak{g}$ into $T$-weight spaces, where $\Phi \subset \mathbb{X}^*(T)$ is the set of roots. For any two $\lambda, \tau \in \mathbb{X}_*(T)_\bQ$, we write
    \[ (\lambda, \tau)_{\mathfrak{g}} := \sum_{\alpha \in \mathbb{X}^*(T)} \dim_k(\mathfrak{g}_{\alpha}) \langle \lambda, \alpha \rangle \langle \tau, \alpha \rangle = \sum_{ \alpha \in \Phi} \langle \lambda, \alpha \rangle \langle \tau, \alpha \rangle.\]
    Note that the bilinear form $(-,-)_{\mathfrak{g}}$ is $W_G$-invariant, symmetric and semidefinite. In fact, the kernel of the bilinear form is the subspace $\mathbb{X}_*(T)_\bQ^{W_G^{\circ}}$ of invariants under the Weyl group $W_G^{\circ} \subset W_G$ of the neutral component $G^{\circ} \subset G$, which is also the subspace $\mathbb{X}_*(Z(G^\circ))_\bQ$ of cocharacters central in $G^\circ$. It follows that the restriction of $(-,-)_{\mathfrak{g}}$ to the kernel $K$ of the averaging operator $\mathbb{X}_*(T)_\bQ \to \mathbb{X}_*(T)_\bQ^{W_G^{\circ}}$ is positive definite. We denote by $\tr_{\mathfrak{g}}: \mathbb{X}_*(T)_\bQ \to \mathbb{X}^*(T)_\bQ$ the corresponding $W_G$-equivariant morphism that sends $\lambda$ to $(\lambda, -)_{\mathfrak{g}}$, so that $\tr_{\mathfrak{g}}(\lambda)=\sum_{\alpha\in \Phi} \langle \lambda,\alpha\rangle \alpha$.
\end{para}

The significance of the trace pairing in our context arises from the following computation. 

\begin{lemma}[Computation of the weight]
    \label{lemma: computation of the weight}
    Let $\lambda\in \mathbb{X}_*(T)$ be a cocharacter and $L=L_G(\lambda)$ its centralizer in $G$. Let $M$ be an $L$-bundle, corresponding to a morphism $p\colon \mathrm{B}\mathbb G_{\mathrm{m}} \to \Bun_G$, and let $[F,d]\in \pi_0(L)_\bQ$ be the rational degree of $M$. We regard $d$ as an element of $\mathbb{X}_*(T)_\bQ$ via the identifications $\mathbb{X}^*(L^F)_\bQ^\vee=\mathbb{X}_*(Z(L^F))_\bQ\subset \mathbb{X}_*(T)_\bQ$. Then $\mathrm{wt}(p^*\mathcal \cL_{\det})=(\lambda, d)_{\mathfrak{g}}$.
\end{lemma}
\begin{proof}
This is a standard Riemann-Roch computation, see for example \cite[1.F.c]{_Heinloth_HilbertMumfordstabilityonalgebraicstacksandapplicationstoGbundlesoncurves} and \cite[Lemma~4.8]{gauged-maps}. The Lie algebra $\mathfrak{g}$ splits as a direct sum $\mathfrak{g}=\bigoplus_{n\in\bZ} \mathfrak{g}_n$, where $\lambda$ acts on $\mathfrak{g}_n$ with weight $n$. Each $\mathfrak{g}_n$ is an $L$-representation. The pullback of the universal bundle $E_{\mathrm{univ}}$ along $p\times \mathrm{id}_C\colon \mathrm{B}\mathbb G_{\mathrm{m}} \times C\to \mathrm{Bun}_C\times C$ is the graded vector bundle $V=\bigoplus_{n\in \bZ} M(\mathfrak{g}_n)$. By flat base change and Riemann-Roch, denoting $q\colon \mathrm{B}\mathbb G_{\mathrm{m}} \times C\to \mathrm{B}\mathbb G_{\mathrm{m}}$ the projection, we have
\begin{align*}
    \mathrm{wt}(p^*\mathcal L_{\det}) &= \mathrm{wt}\left(\det\left(Rq_*V\right)\right) = \sum_{n\in \bZ} n \cdot \chi(M(\mathfrak{g}_n)) \\
    &=\sum_{n\in \bZ} n \cdot \left(\deg(M(\mathfrak{g}_n)) + \mathrm{rk}(M(\mathfrak{g}_n))(1-g)\right) = \sum_{n\in \bZ} n \cdot \deg(M(\mathfrak{g}_n))=(\lambda, d)_{\mathfrak{g}},
\end{align*}
as desired.
\end{proof}

\begin{para}[Levi subgroups] \label{para: parabolic subgroups}
    Given a parabolic subgroup $P \subset G$ containing $T$, let $L$ denote its corresponding Levi quotient. We view $L$ as a subgroup of $P$ by considering the unique Levi subgroup $L \subset P$ which contains the maximal torus $T$. The inclusion $L \hookrightarrow G$ induces a natural inclusion $W_L \hookrightarrow W_G$. 

\end{para}

\begin{para}[Dominant cocharacters and characters] \label{para: dominant cocharacters}
    Given a parabolic subgroup $P \subset G$ containing the maximal torus $T$, we say that a rational cocharacter $\lambda\in \mathbb{X}_*(T)_\bQ$ is $P$-dominant if $P_G(\lambda)= P$. We say that a rational character $\mathbb{X}^*(T)_\bQ$ is $P$-dominant if it is of the form $\tr_{\mathfrak{g}}(\lambda)$ for some $P$-dominant rational cocharacter $\lambda \in \mathbb{X}_*(T)_\bQ$. By construction, $P$-dominant rational cocharacters lie in the subspace of invariants $\mathbb{X}_*(T)_\bQ^{W_L}$, where $L \subset P$ is the Levi subgroup as in \cref{para: parabolic subgroups}. Similarly, dominant rational characters lie on $\mathbb{X}^*(T)_\bQ^{W_L}$, and hence we may view them canonically as elements of $\mathbb{X}^*(L)_\bQ = \mathbb{X}^*(P)_\bQ$ by \Cref{theorem: rational characters and Weyl group}.
\end{para}

\begin{para}[Degree of parabolic bundles]
    Let $P \subset G$ be a parabolic subgroup containing $T$, with corresponding Levi subgroup $L \subset P$ as in \cref{para: parabolic subgroups}. Given a $P$-bundle $M$, we denote by $[F_{M}, d_{M}] \in \pi_1(L)_{\mathbb{Q}}$ the rational degree of the corresponding associated $L$-bundle $M(L)$ (obtained via the quotient morphism $P \to L$). We may view $d_{M}$ as an element of $\mathbb{X}_*(T)_\bQ$ via the canonical identifications $\mathbb{X}^*(L')_\bQ^{\vee} = \mathbb{X}_*(Z(L'))_\bQ \subset \mathbb{X}_*(T)_\bQ$, where $L'=L^{F_M}$.
\end{para}

\begin{para}[Semistability in terms of parabolic reductions] \label{para: semistability in terms of parabolic reductions}
    Using the computation of the weight of $\cL_{\det}$ in \cref{lemma: computation of the weight} and the description of filtrations in \cref{para: filtrations and graded points}, one can see that our notion of semistability in \Cref{defn: semistability bundles} can be expressed in terms of parabolic reductions as follows. Let $K \supset k$ be a field extension, and let $E \in \Bun_G(K)$ be a $G$-bundle on $C_K$. Then, $E$ is semistable if and only if for all parabolic subgroups $P \subset G$ containing $T$, all $P$-reductions of structure group $E_P$ of $E$ and all $P$-dominant rational cocharacters $\lambda \in \mathbb{X}_*(T)_\bQ$, we have $(\lambda, d_{E_P})_{\mathfrak{g}} \leq 0$.

    We may equivalently express this in terms of dominant characters: $E$ is semistable if and only if for all parabolic subgroups $P \subset G$ containing $T$, all $P$-reductions of structure group $E_P \subset E$ and all $P$-dominant rational characters $\chi$, we have $\deg(E_P(\chi)) \leq 0$, where we view $\chi \in \mathbb{X}^*(T)_\bQ^{W_{L}}$ as a rational character of $P$ as in \cref{para: dominant cocharacters}.
\end{para}

\subsection{Semistability under change of group}


\begin{theorem}[Semistability under change of group] \label{thm: semistability under change of group} Let $f: G \to H$ be a homomorphism of smooth linearly reductive affine algebraic groups over $k$. Let $K \supset k$ be a field extension, and let $E \in \Bun_G(K)$ be a $G$-bundle on $C_K$.
\begin{enumerate}
    \item If the rational degree $[F, d_E]$ of $E$ is not $f$-adapted, then the associated $H$-bundle $E(H)$ is not semistable.
    \item If the rational degree $[F, d_E]$ of $E$ is $f$-adapted (as in \Cref{definition: adapted degree}) and $E$ is semistable, then the associated $H$-bundle $E(H)$ is semistable.
\end{enumerate}
\end{theorem}

We provide the proof of \Cref{thm: semistability under change of group} below in \cref{proof of semistability under change of group}. Before doing that, let us record a remark, a corollary, and some examples of the theorem.

\begin{remark}[Comparison with earlier results]
    To our knowledge, the statement of \Cref{thm: semistability under change of group} above is new even in the case when $G$ and $H$ are connected (see \Cref{example: the case of connected groups}). It gives a complete characterization of when the associated bundle construction preserves the semistability of a principal bundle. The result that is usually stated in the literature for connected reductive groups, going back to Ramanathan's thesis \cite[Prop. 3.17]{ramanathan_moduli_principal_bundles} and the article by Ramanan and Ramanathan \cite[Thm. 3.18]{ramanan-ramanathan}, is the forthcoming weaker \Cref{coroll: semistability under change of group weaker}. 
\end{remark}

\begin{corollary} \label{coroll: semistability under change of group weaker}
    Let $f: G \to H$ be a homomorphism of smooth linearly reductive affine algebraic groups over $k$. Assume that the image of the maximal central torus $Z_{G^{\circ}}^{\circ}$ of the neutral component $G^{\circ} \subset H$ lies inside the maximal central torus $Z_{H^{\circ}}^{\circ}$ of the neutral component $H^{\circ} \subset H$. Let $K \supset k$ be a field extension, and let $E \in \Bun_G(K)$ be a semistable $G$-bundle on $C_K$. Then, the associated $H$-bundle $E(H)$ is semistable.
\end{corollary}
\begin{proof}
    This is an immediate consequence of \Cref{thm: semistability under change of group}, because the assumption $f(Z_{G^{\circ}}^{\circ}) \subset Z_{H^{\circ}}^{\circ}$ ensures that every rational degree $[F,d] \in \pi_1(G)_{\mathbb{Q}}$ is $f$-adapted.
\end{proof}

\begin{example}[The case of connected groups] \label{example: the case of connected groups}
    For ease of reference, let us spell out \Cref{thm: semistability under change of group} in the case when both groups $G$ and $H$ are connected, as in \Cref{example: rational degrees for connected groups}. In this case, given a bundle $E$, the rational degree $d_E$ is the linear functional in $\pi(G)_{\mathbb{Q}} = \mathbb{X}_*(T)_\bQ^{W_G}$ described in \Cref{example: rational degrees for connected groups}. Consider the morphism $\tau: \pi_1(G)_{\mathbb{Q}} = \mathbb{X}_*(T)_\bQ^{W_G} \hookrightarrow \mathbb{X}_*(T)_\bQ \to \mathbb{X}_*(T')_\bQ$. Then, \Cref{thm: semistability under change of group} states the following complete characterization of when $f_*: \Bun_G \to \Bun_H$ preserves semistability:
    \begin{enumerate}
        \item Suppose that $\tau(d_E) \notin \mathbb{X}_*(T')_\bQ^{W_H}$ (i.e. $d_E$ is not $f$-adapted). Then, the image of the morphism $f_*: \Bun_G^{d_E} \to \Bun_H$ lies in the unstable locus of $\Bun_H$.

        \item Otherwise, suppose that $\tau(d_E) \in \mathbb{X}_*(T')_\bQ^{W_H}$. Then the morphism $f_*: \Bun_G^{d_E} \to \Bun_H$ preserves semistability; it sends semistable $G$-bundles to semistable $H$-bundles.

    \end{enumerate}
\end{example}

\begin{example}[Direct sum homomomorphism] \label{example: semistability under direct sum homomorphism}
    In the context of \Cref{example: rational degrees for direct sum homomorphism}, \Cref{thm: semistability under change of group} reduces to a familiar statement, namely:
    \begin{enumerate}
        \item Let $(E_1, E_2)$ be a pair of vector bundles. If the slopes of $E_1$ and $E_2$ are not equal, then $E_1 \oplus E_2$ is not semistable.
        \item Let $(E_1, E_2)$ be a pair of semistable vector bundles. If the slopes of $E_1$ and $E_2$ are equal, then the direct sum $E_1 \oplus E_2$ is semistable.
    \end{enumerate}
\end{example}

\begin{example}\label{example: subgroups of pi0 are necessary}
Let $G=N_{\GL_2}(T)$ be the normalizer of the standard maximal torus $T=\G_{m,k}^2$ inside $\GL_2$. We have $\pi_0(G) = S_2 = \mathbb{Z}/2\mathbb{Z}$, which acts by conjugation on $\mathbb{X}_*(T)_\bQ= \mathbb{Q}^{\oplus 2}$ via the standard permutation representation (switching the coordinates). We have an identification
\[\pi(G)_\Q=(\Q^{\oplus 2}/S_2) \sqcup \Q,\]
where the set of $S_2$-orbits $\Q^{\oplus 2}/S_2=\mathbb{X}_*(T)_\bQ/S_2$ corresponds to the trivial subgroup $F=\{1\}\subset S_2$, and $\Q=\mathbb{X}^*(G)_\bQ^\vee = \mathbb{X}_*(T)_\bQ^{S_2}$ corresponds to the whole subgroup $F=S_2\subset S_2$.
Consider an $G$-bundle $E$ that is induced via $T\to G$ from a pair $(L_1,L_2)$ of line bundles on $C$, and let $d_1,d_2$ be the degrees of $L_1$ and $L_2$. The rational degree of $E$ is $[(d_1,d_2)]\in \Q^{\oplus 2}/S_2$, and it is adapted for the inclusion $G\to \GL_2$ precisely if $d_1=d_2$.

If we had defined the notion of rational degree without taking into account subgroups of $\pi_0(G)$ (just as a functional in $\mathbb{X}^*(G)_\bQ^{\vee}$), then the rational degree of $E$ would just be $d_1+d_2\in \mathbb{X}^*(G)_\bQ^\vee$, which does not contain enough information to determine whether $d_1=d_2$, that is, whether the associated bundle $E(\GL_2)$ is semistable.
\end{example}

\begin{para} \label{proof of semistability under change of group}
    \begin{proof}[Proof of \Cref{thm: semistability under change of group}] For this proof, we denote by $\mathfrak{h}$ the Lie algebra of $H$. We fix a maximal torus $T' \subset H$ containing $f(T)$. We denote by $\tau: \mathbb{X}_*(T)_\bQ \to \mathbb{X}_*(T')_\bQ$ the induced morphism of vector spaces.
    
    Let us make some preliminary reductions and observations before proceeding with the proof. We may choose a lift $\widetilde{E}$ to $\Bun_{G^{F}}$ and replace $G$ and $H$ with $G^{F}$ and $H^{f(F)}$ respectively. Hence, in view of \Cref{prop: behavior stability under change of group of components}, we may assume without loss of generality that $F = \pi_0(G)$ and $f(\pi_0(G)) = \pi_0(H)$ throughout this proof. In this case, there is a unique representative $(F,d_E)$ for $[F,d_E]$.
    
    If $P \subset G$ is a parabolic subgroup containing $T$ such that $E$ admits a $P$-reduction of structure group $E_P$, then notice that there is an inclusion $\pi_0(P) \hookrightarrow \pi_0(G)$, and that $E_P(\pi_0(P))$ is a $\pi_0(P)$-reduction of structure group of $E(\pi_0(G))$. Since we are assuming that $E(\pi_0(G))$ is connected, it follows that $\pi_0(P) = \pi_0(L) = \pi_0(G)$, and that the corresponding $\pi_0(P)$-bundle $E_P(\pi_0(P)) = E_L(\pi_0(L))$ is connected. Hence, we get again that there is a unique representative $d_{E_P}$ for the degree of $E_P$. A similar reasoning shows that for any parabolic subgroup $Q \subset H$ containing $T'$ such that there is a $Q$-reduction $E(H)_Q$ of $E(H)$, we have an equality $\pi_0(Q)= \pi_0(H)$ and there is a unique representative $d_{E_Q(H)}$ for the rational degree. For the rest of the proof, we use freely  the uniqueness of representatives for rational degrees of parabolic reductions.
    
    \noindent\emph{Proof of part (i).} Suppose that $[F_E, d_E] = [\pi_0(G), d_E]$ is not $f$-adapted. Set $\lambda=d_E$, seen as an element of $\mathbb{X}_*(T)_\bQ$. Then, the image $\lambda'=\tau(\lambda)$ of $\lambda$ under the morphism 
    \[\tau: \mathbb{X}^*(G)_\bQ^{\vee} = \mathbb{X}_*(T)_\bQ^{W_{G}} \to \mathbb{X}_*(T')_\bQ\] 
    does not lie in the kernel $\mathbb{X}_*(T')_\bQ^{W_{H^{\circ}}}$ of the semidefinite bilinear form $(-,-)_{\mathfrak{h}}$ on $\mathbb{X}^*(T')_\bQ$. Let $\lambda' = \tau(\lambda)$ denote the corresponding rational cocharacter of $T'$, and let $Q\coloneq P_H(\lambda')$ be the corresponding parabolic subgroup of $H$, with unique Levi subgroup $L' = L_H(\lambda')$ containing $T'$. Notice that we have $L_G(\lambda) = P_G(\lambda) = G$, because $\lambda=d_E \in \mathbb{X}_*(T)_\bQ^{W_G}$ is central (\Cref{theorem: cocharacters and Weyl group}). Hence, get a morphism $G= L_G(\lambda) \to L_H(\lambda') = L' \hookrightarrow Q$. The associated $Q$-bundle $E(Q)$ is a $Q$-reduction of structure group of $E(H)$. Note that $\lambda'$ is invariant under the Weyl group $W_{L'}$ of $L'$, since being is central in $L'$. By the description of the map $\mathbb{X}^*(G)_\bQ^\vee\to \mathbb{X}^*(L')_\bQ^\vee$ in \cref{para: decomposition of map on rational degrees}, it follows that $\lambda'$ equals the degree $d_{E(Q)}$ of $E(Q)$ as elements of $\mathbb{X}_*(T')_\bQ$. We have that $\lambda'$ is $Q$-dominant by construction, and
    \[ (\lambda', d_{E(Q)})_{\mathfrak{h}} = (\lambda', \lambda')_{\mathfrak{h}} >0.\]
    Hence, the data of parabolic group $Q$, parabolic reduction $E(Q)$ and dominant rational cocharacter $\lambda'$ witness the instability of $E(H)$ as in \cref{para: semistability in terms of parabolic reductions}.

\noindent\emph{Proof of part (ii).} We use the Ramanan--Ramanthan construction \cite{ramanan-ramanathan}, see also \cite[Construction 2.33]{herrero2023meromorphic}. Suppose for the sake of contradiction that the associated bundle $E(H)$ is unstable. This means that there is a parabolic subgroup $Q \subset H$, a parabolic reduction $E(H)_Q \subset E$ and a $Q$-dominant cocharacter $\lambda' \in \mathbb{X}_*(T')_\bQ$ such that $(\lambda', d_{E(H)_Q})_{\mathfrak{h}} >0$.

The reduction of structure group $E(H)_Q$ corresponds to a section $\sigma: C \to E(H/Q)$ of the associated $H/Q$-fiber bundle $E(H/Q)$. The morphism $G \to H$ induces an action of $G$ on the flag scheme $H/Q$. We have a $Q$-dominant rational character $\chi:= \tr_{\mathfrak{h}}(\lambda') \in\mathbb{X}^*(T')_\bQ^{W_{L'}} = \mathbb{X}^*(Q)_\bQ$, where $L'$ is the unique Levi subgroup of $Q$ containing $T'$. This induces a $G$-equivariant ample line bundle $\cO(-\chi)$ on $H/Q$. Consider the kernel $K \subset \mathbb{X}_*(T)_\bQ$ of the averaging operator $\mathbb{X}_*(T)_\bQ \twoheadrightarrow \mathbb{X}_*(T)_\bQ^{W_{G^{\circ}}}$. We have $\mathbb{X}_*(T)_\bQ = \mathbb{X}_*(T)_\bQ^{W_{G^{\circ}}} \oplus K$, and $\mathbb{X}_*(T)_\bQ^{W_{G^{\circ}}}$ is the kernel of the semidefinite bilinear form $(-,-)_{\mathfrak{g}}$. Similarly as in \cite[Construction 2.33]{herrero2023meromorphic}, we fix a choice of $W_G$-invariant positive definite symmetric bilinear form $b'$ on $\mathbb{X}_*(T)_\bQ^{W_{G^{\circ}}}$ and consider the corresponding $W_G$-invariant positive definite bilinear form $b= b' \oplus (-,-)_{\mathfrak{g}}|_K$ on $\mathbb{X}_*(T)_\bQ = \mathbb{X}_*(T)_\bQ^{W_{G^{\circ}}} \oplus K$ (the choice of $b'$ will ultimately not matter).

The construction in \cite[Construction 2.33]{herrero2023meromorphic} applies verbatim thanks to \Cref{thm: description of centres} to yield a parabolic subgroup $P\subset G$ with Levi subgroup $T \subset L \subset P$, a parabolic reduction $E_P \subset E$ and a $P$-dominant ``maximally destabilizing cocharacter'' $\lambda \in \mathbb{X}_*(T)$ such that $\xi := (\lambda,-)_b \in \mathbb{X}^*(T)_\bQ^{W_L} = \mathbb{X}^*(P)_\bQ$ satisfies $\deg(E_P(\xi)) > 0$.

Let $\kappa \subset \mathbb{X}_*(T)_\bQ$ denote the subspace of elements $\sigma$ such that $\tau(\sigma)$ lies in $\mathbb{X}_*(T')_\bQ^{W_{H^{\circ}}}$. We note that, since $\pi_0(Q) = \pi_0(H)$, the inclusion $H^{\circ}/Q^{\circ} \hookrightarrow H/Q$ is an isomorphism. It follows that the centre $Z(H^{\circ})$ acts trivially on $H/Q$. By \Cref{theorem: cocharacters and Weyl group}, for all $\sigma \in \kappa$ the condition $\tau(\sigma) \in \mathbb{X}_*(T')_\bQ^{W_{H^{\circ}}}$ ensures that $\sigma$ acts trivially on $H/Q$ (when viewed as cocharacter of $G$ up to a power). On the other hand, since $\chi = \tr_{\mathfrak{h}}(\lambda)$ pairs trivially with every element in the subspace $\mathbb{X}^*(T')_\bQ^{W_{H^{\circ}}}$, we have $\langle \tau(\sigma), \chi \rangle =0$ for all $\sigma \in \kappa$.
 
 Since the rational cocharacters in $\kappa$ act trivially on $H/Q$ and have zero pairing with $\chi$, the fact that $\lambda$ is maximally destabilizing implies that $\lambda$ is in the orthogonal complement $\kappa^{\perp}$ of $\kappa$ with respect to the bilinear form $b$. Our assumption that $[\pi_0(G), d_E]$ is $f$-adapted means that $d_E \in \kappa$, and hence $\lambda \in (\mathbb{Q} d_E)^{\perp}$. 

Let us consider the first projection $\pr_1(\lambda) \in \mathbb{X}_*(T)_\bQ^{W_{G^{\circ}}}$ of $\lambda \in \mathbb{X}_*(T)_\bQ = \mathbb{X}_*(T)_\bQ^{W_{G^{\circ}}} \oplus K$. Recall that the existence of the $P$-reduction $E_P$ of $E$ implies that $\pi_0(P) = \pi_0(L) = \pi_0(G)$. This implies that $W_L$ and $W_{G^{\circ}}$ jointly generate $W_G$. Since the cocharacter $\pr_1(\lambda)$ is invariant under $W_L$ (because $\lambda$ is so) and under $W_{G^{\circ}}$, we conclude that it is invariant under the full Weyl group $W_G$. In particular, $\lambda - \pr_1(\lambda)$ is also $P$-dominant, since subtracting $W_G$-invariant elements has no effect on $P$-dominance (see  \Cref{theorem: cocharacters and Weyl group}). We may write $\xi = \xi_1 + \xi_2$, where we set: 
\[\xi_1 = (\pr_1(\lambda),-)_b, \text{  and    } \; \; \; \xi_2 = (\lambda - \pr_1(\lambda),-)_b = (\lambda - \pr_1(\lambda),-)_{\mathfrak{g}},\]
where the last equality follows because of the agreement of the restrictions of $b$ and $(-,-)_{\mathfrak{g}}$ to $K$ and the fact that $\lambda - \pr_1(\lambda) \in K$. Since $\lambda - \pr_1(\lambda)$ is $P$-dominant, it follows from the last displayed equality above that $\xi_2$ is a $P$-dominant rational character. To conclude the proof, we shall show that $\deg(E_P(\xi_2))>0$; this would contradict the semistability of $E$ in view of \cref{para: semistability in terms of parabolic reductions}.

The equality $\deg(E_P(\xi))>0$ implies that $\deg(E_P(\xi_1)) + \deg(E_P(\xi_2))>0$. It suffices to show then that $\deg(E_P(\xi_1)) = 0$. Unraveling the definitions, we see that 
\[\deg(E_P(\xi_1)) = \langle d_{E_P}, \xi_1\rangle = ( \pr_1(\lambda), d_{E_P})_b.\]
Since $b$ and $\pr_1(\lambda)$ are $W_G$-invariant, we see that $( \pr_1(\lambda), d_{E_P})_b=( \pr_1(\lambda), d_{E_P}^{W_G})_b$, where we denote by $d_{E_P}^{W_G}$ the projection of $d_{E_P} \in \mathbb{X}_*(T)_\bQ$ under the averaging operator $\mathbb{X}_*(T)_\bQ \to \mathbb{X}_*(T)_\bQ^{W_G}$. It follows from the definition of degree that the projection $d_{E_P}^{W_G}$ is equal to $d_E$. Hence we get
\[ \deg(E_P(\xi_1)) = ( \pr_1(\lambda), d_{E})_b.\] 
Since $\lambda$ lies in $(\mathbb{Q}d_E)^{\perp}$ and the direct summands in the decomposition $\mathbb{X}^*(T)_\bQ = \mathbb{X}^*(T)_\bQ^{W_{G^{\circ}}} \oplus K$ are $b$-orthogonal by construction, it follows that the projection $\pr_1(\lambda)$ is also in $(\mathbb{Q}d_E)^{\perp}$. We conclude that $( \pr_1(\lambda), d_{E})_b=0$, which implies that $\deg(E_P(\xi_1))=0$, as desired.
\end{proof} 
\end{para}

\begin{proposition}[Semistability and adjoint bundles] \label{prop: equivalence semistability notions}
    Let $\mathfrak{g}$ denote the Lie algebra of $G$. Let $K \supset k$ be a field extension, and let $E \in \Bun_G(K)$ be a $G$-bundle on $C_K$. Then, $E$ is semistable if and only if the associated adjoint vector bundle $E(\mathfrak{g})$ is semistable.
\end{proposition}
\begin{proof}
    \Cref{coroll: semistability under change of group weaker} immediately implies one direction: if $E$ is semistable, then the associated vector bundle $E(\mathfrak{g})$ is semistable. 
    
    The proof of the converse is simpler: one may argue similarly as in \cite[Corollary 2.32]{herrero2023meromorphic} using our alternative description of semistability in terms of parabolic reductions in \cref{para: semistability in terms of parabolic reductions}; the key is to notice that the trace pairing $\tr_{\mathfrak{gl}(\mathfrak{g})}$ of the adjoint representation $\mathfrak{gl}(\mathfrak{g})$ of $\GL(\mathfrak{g})$ agrees up to a positive scaling with the trace pairing $\tr_{\mathfrak{g}}$ of the standard representation $\mathfrak{g}$ of $\GL(\mathfrak{g})$ when restricted to the kernel of the projection of the averaging operator by $W_{\GL(\mathfrak{g})}$.
\end{proof}

\begin{corollary}
    Let $f\colon L=L_G(\lambda)\to G$ be the inclusion of a Levi subgroup for a one-parameter subgroup $\lambda$ of $G$, and let $E$ be an $L$-bundle. Then $E(G)$ is semistable if and only if the degree $[F_E,d_E]$ of $E$ is $f$-adapted and $E$ is semistable.
\end{corollary}

\begin{proof}
    If $[F_E,d_E]$ is $f$-adapted and $E$ is semistable, then the associated bundle $E(G)$ is semistable by \Cref{thm: semistability under change of group}.

    Conversely, assume that $E(G)$ is semistable. Then $[F_E,d_E]$ is $f$-adapted by \Cref{thm: semistability under change of group}. By linear reductivity of $L$, the Lie algebra $\mathfrak{l} \coloneq \mathrm{Lie}(L)$ is a direct summand of $\mathfrak{g}$ as representations of $L$ under the adjoint action. This implies that the adjoint bundle $E(\mathfrak{l})$ is a direct summand of $E(\mathfrak{g}) = E(G)(\mathfrak{g})$. We conclude by \cref{coroll: semistability under change of group weaker} and the fact that a direct summand of a semistable vector bundle is semistable.
\end{proof}

\subsection{Comparison of existing notions of semistability} \label{section: comparison notions semistability}

In this subsection we survey the existing notions of semistability for principal bundles in the literature. We show that they all agree, thus answering positively a question by \textcite{biswas-gomez-kahler-einstein}.

\begin{para}[Hilbert-Mumford semistability]
    The notion of Hilbert-Mumford semistability in \Cref{defn: semistability bundles} agrees with the one in \textcite{moduli-principal-bundles-large}, which is used to construct a projective moduli space of semistable principal bundles via Geometric Invariant Theory. As explained in \cref{para: semistability in terms of parabolic reductions}, there is an alternative description Hilbert-Mumford semistability in terms of parabolic reductions. The description of semistability in \cref{para: semistability in terms of parabolic reductions} agrees, in the case of connected reductive groups, with the classical notion of semistability defined by \textcite{ramanathan_moduli_principal_bundles}; this was first observed by \textcite{_Heinloth_HilbertMumfordstabilityonalgebraicstacksandapplicationstoGbundlesoncurves}.
\end{para}

\begin{para}[ad-semistability]
    We say that a $G$-bundle $E$ is ad-semistable if and only if the associated adjoint vector bundle $E(\mathfrak{g})$ is slope semistable. In the case of connected reductive groups, ad-semistability was considered by Atiyah-Bott in \cite[\S 10]{_Atiyah_TheYangMillsEquationsoverRiemannSurfaces}. We note that for connected reductive groups it was shown that ad-semistability coincides with the open stratum of the stratification induced by the Yang-Mills functional, see \cite[Prop. 10.9]{_Atiyah_TheYangMillsEquationsoverRiemannSurfaces}. In the case of disconnected reductive groups, ad-semistability is the notion considered in \textcite{moduli-disconnected-groups} to build a moduli space of principal bundles.
\end{para}

\begin{para}[Einstein-Hermite semistability]
    \textcite[Definition 1]{biswas-gomez-kahler-einstein} define another notion of semistability, which for the purposes of this paper we shall call Einstein-Hermite semistability. Given a $G$-bundle $E$, the quotient $Y:= E/G^{\circ} \to C$ by the neutral component $G^{\circ} \subset G$ is a finite \'etale cover of $C$. By construction, $Y$ classifies reductions of structure group of $E$ to $G^{\circ}$, and hence the pullback $E|_{Y}$ admits a canonical reduction of structure group $E_{G^{\circ}}$. The principal bundle $E$ is called Einstein-Hermite semistable if for all connected components $Y' \subset Y$ the $G^{\circ}$-bundle $E_{G^{\circ}}|_{Y'}$ is semistable in the sense of Ramanathan (or, equivalently, is Hilbert-Mumford semistable). As usual, one may define a corresponding notion of Einstein-Hermite polystability, see . It is shown in \cite{biswas-gomez-kahler-einstein} that the principal bundle $E$ is Einstein-Hermite polystable if and only if it admits an Einstein-Hermitian structure.
\end{para}

It is left as an open question in \cite{biswas-gomez-kahler-einstein} whether Einstein-Hermite semistability agrees with Hilbert-Mumford semistability. We answer this positively in the following.

\begin{proposition} \label{prop: survey all semistability conditions}
    For any given $G$-bundle $E$ on the curve $C$, the following are equivalent:
    \begin{enumerate}
        \item $E$ is Hilbert-Mumford semistable.
        \item $E$ is ad-semistable.
        \item $E$ is Einstein-Hermite semistable.
    \end{enumerate}
\end{proposition}
\begin{proof}
    The equivalence of (ii) and (iii) is proven in \cite[Lemma 2]{biswas-gomez-kahler-einstein}. The equivalence of (i) and (ii) is the content of \Cref{prop: equivalence semistability notions}.
\end{proof}

\preparebibliography
\printbibliography

\authorinforule

\authorinfo{Andres Fernandez Herrero}{andresfh@sas.upenn.edu}{Department of Mathematics, University of Pennsylvania,
Philadelphia, PA 19104,USA}

\authorinfo{Andrés Ibáñez Núñez}
    {andres.ibaneznunez@columbia.edu}
    {Mathematics Department, Columbia University, New York, NY 10027, USA}
    
\end{document}